\documentclass{article}

\usepackage{PRIMEarxiv}

\usepackage[utf8]{inputenc} 
\usepackage[T1]{fontenc}    
\usepackage{hyperref}       
\usepackage{url}            
\usepackage{booktabs}       
\usepackage{amsfonts}       
\usepackage{nicefrac}       
\usepackage{microtype}      
\usepackage{lipsum}
\usepackage{fancyhdr}       
\usepackage{graphicx}       
\graphicspath{{media/}}     
\usepackage{siunitx}
\usepackage{amsfonts}
\usepackage{amsmath}
\usepackage{amsthm}
\usepackage{enumitem} 
\usepackage{url}

\newtheorem{theorem}{Theorem}

\newtheorem{definition}{Definition}
\newtheorem{example}{Example}

\newtheorem{remark}{Remark}

\title{On the    Generalized  Hukuhara Nabla Differentiability of Fuzzy Functions  on Time Scales via Characterization Theorem
\thanks{\textit{\underline{Citation}}: 
\textbf{Authors. Title. Pages.... DOI:000000/11111.}} 
}

\author{
  Funda Raziye MERT \\
  Department of Software Engineering \\
  Adana Alparslan Türkeş Science and Technology University \\
  Adana, TÜRKİYE\\
  \texttt{rmert@atu.edu.tr} \\
   \And
  Selami BAYEĞ \\
  Department of Industrial Engineering\\
  University of Turkish Aeronautical Association \\
  Ankara, TÜRKİYE\\
  \texttt{sbayeg@thk.edu.tr} \\
   \AND
   Billur KAYMAKÇALAN \\
   Department of Computer Engineering\\
   University of Turkish Aeronautical Association \\
  Ankara, TÜRKİYE\\
  \texttt{bkaymakcalan@thk.edu.tr} \\
}

\begin{document}
\maketitle

\begin{abstract}
In this paper, we provide a comprehensive investigation of the  generalized Hukuhara nabla derivative for fuzzy functions on time scales. We establish some characterizations of  generalized Hukuhara  nabla differentiable fuzzy functions on time scales via the nabla differentiability of their endpoint functions. These characterizations complete and extend previous findings of fuzzy number valued functions from the real domain to those on time scales.  Additionally, we generalize the product rule, previously proven for interval valued functions in the real case, to fuzzy number valued functions on time scales.
\end{abstract}

\keywords{Time scale \and Generalized Hukuhara  nabla derivative \and  Fuzzy functions \and  Generalized Hukuhara difference }

\section{Introduction}
In recent years, the study of dynamic systems has increasingly required mathematical tools capable of handling both continuous and discrete phenomena. Traditional calculus methods, such as differential and difference equations, typically operate within distinct frameworks: one for continuous systems and another for discrete ones. However, many real-world systems exhibit hybrid behavior, transitioning between continuous and discrete states. To address this challenge, the theory of time scales, introduced by Stefan Hilger in the 1980s, provides a unified approach for analyzing dynamic systems on arbitrary time scales 
 \cite{Hilger1}.
\\ \\
Nabla  time scale calculus   is  a central component of  time scale theory.  Specifically,  it deals with backward differences and provides a powerful tool for investigating problems in a variety of fields, including biology, economics, engineering, and physics. The development of nabla time scale calculus, particularly the nabla derivative and integral, enables a more generalized analysis of dynamic equations, offering significant flexibility in modeling complex systems \cite{Peterson, Peterson1, Guseinov, Anderson}.
\\ \\
Fuzzy set theory, introduced by Lotfi Zadeh in 1965, extends classical set theory by allowing elements to have degrees of membership rather than strict binary classification \cite{Zadeh}. Unlike traditional sets, where an element either belongs or does not belong, fuzzy sets assign a membership value between 0 and 1, representing the degree of belonging. This flexibility makes fuzzy set theory particularly useful for handling uncertainty, vagueness, and imprecise information in various fields, including artificial intelligence, control systems, decision-making, and pattern recognition. By modeling real-world problems with gradual transitions rather than rigid boundaries, fuzzy set theory provides a powerful tool for dealing with complex and ambiguous data \cite{Klir, Gomes, Bedebook}.
\\ \\
The generalized Hukuhara difference is an extension of the traditional Hukuhara difference, used to analyze the difference between fuzzy numbers in fuzzy set theory. This generalized approach provides a more versatile framework for operations on fuzzy numbers, especially in the context of dynamic systems and differential equations. It allows for a deeper exploration of fuzzy systems across both continuous and discrete time domains. Specifically, the investigation of the generalized Hukuhara nabla derivative has developed into a significant extension of the standard nabla derivative. By combining the generalized  Hukuhara nabla derivative with nabla calculus, this framework provides a powerful tool for investigating systems where traditional difference or differential operators are inadequate, such as those involving fuzzy processes or hybrid time scales. For recent studies on nabla Hukuhara derivative for fuzzy functions on time scales, we refer to the references
\cite{Leelavathi2019, LKM1, Leelavathi2020, Leelavathi2018,Leelavathi2019-3, Leelavathi2019-4, you, Hari}. In \cite{Leelavathi2019, LKM1}, Leelavathi et al. introduced the nabla Hukuhara derivative and the second-type nabla Hukuhara derivative for fuzzy functions on time scales under the Hukuhara difference. They proved the existence and uniqueness of these derivatives and established their fundamental properties. Additionally, they examined the nabla Hukuhara and the second type nabla Hukuhara derivatives of scalar multiplication, sum, and product of two fuzzy functions on time scales. In \cite{Fard}, Fard et al. introduced the $\Delta$-Hukuhara derivative and the $HK_F$-$\Delta$-integral for fuzzy functions on time scales and proved some of their basic properties. In a manner similar to \cite{Fard}, You et al. introduced the nabla Hukuhara derivative of fuzzy functions on time scales and showed some of its basic properties in \cite{you}. In \cite{Hari}, Hari et al. proved the local existence and uniqueness theorem for fuzzy nabla-dynamic equations on time scales under the characterization theorem by utilizing the strongly generalized nabla Hukuhara differentiability to 
transform FNDES-$\tau$ into a system of nabla-dynamic equations on time scales (NDES-$\tau$). In \cite{Cano}, Chalco-Cano et al. discussed several characterizations of generalized Hukuhara differentiability for fuzzy functions, relying on the differentiability of their endpoint functions. However, in \cite{qiu}, Qiu presented a counterexample demonstrating that the generalized Hukuhara differentiability of an interval valued function at a given point does not always align with the one-sided differentiability of its endpoints. Following this, he developed a more complete framework for characterizing generalized Hukuhara differentiability in such cases. Later, in \cite{longo}, Longo et al. expanded these characterization results to cover fuzzy number valued functions.
\\ \\
In this paper, we comprehensively investigate the generalized Hukuhara nabla derivative, a concept introduced by You et al. \cite{you}, with the aim of addressing unresolved gaps in the existing literature while significantly expanding and refining the scope of prior results. We provide several characterizations of generalized  Hukuhara nabla  differentiable fuzzy functions on time scales, based on the nabla differentiability of their endpoint functions. These characterizations serve as valuable tools for the calculus of derivatives of fuzzy functions on time scales. Our results complete and  extend the findings in \cite{Cano} from fuzzy number valued functions in the real domain to those on time scales. We have also extended the product rule, which was proven in \cite{TZ} for interval valued functions in the real case to fuzzy number valued functions on time scales.


\section{Preliminaries}
The following definitions are presented in Bohner et al.  \cite{Peterson} and Agarwal et al. \cite{Agarwal}. A time scale $\mathbb{T}$ is defined as any closed subset of $\mathbb{R}$. The forward and backward jump operators $\sigma, \rho: \mathbb{T}\to\mathbb{T}$  are defined as follows:
$$\sigma(t)=\inf\{s\in\mathbb{T}: s > t\}\,\,\,\text{and}\,\,\,\rho(t)=\sup\{s\in\mathbb{T} : s < t\}$$
(supplemented by $\inf\emptyset:=\sup\mathbb{T}$ and $\sup\emptyset:=\inf\mathbb{T}$). A point $t\in\mathbb{T}$ is
called left-dense if $t>\inf\mathbb{T}$ and $\rho(t)=t$, left-scattered if $\rho(t)<t$, right-dense if $t<\sup\mathbb{T}$
and $\sigma(t)=t$, right-scattered if $\sigma(t)>t$. If $\mathbb{T}$ has a right-scattered minimum $m$, then define
$\mathbb{T}_{\kappa}:=\mathbb{T}-\{m\}$; otherwise, set $\mathbb{T}_{\kappa}:=\mathbb{T}$. The backward graininess $\nu:\mathbb{T}_{\kappa}\to[0,\infty)$
is defined by
$$\nu(t)=t-\rho(t)$$
and it represents the difference between
$t$ and its backward jump 
$\rho(t).$
For $f:\mathbb{T}\to\mathbb{R}$ and $t\in\mathbb{T}_{\kappa}$,  the nabla derivative  of $f$ at $t$, denoted 
$\nabla f(t)$,  is defined as the number (provided it exists) with the property that for any $\varepsilon>0$, there exists a neighborhood $U$ of $t$ such that
$$|[f(\rho(t))-f(s)]-\nabla f(t)[\rho(t) - s]| \leq \epsilon |\rho(t)-s|$$
for all \( s \in U \) \cite{Guseinov}.
\begin{theorem}\label{nabladerrule}\cite{Peterson, Guseinov}
Let \( f: \mathbb{T} \to \mathbb{R} \) be a function, and let \( t \in \mathbb{T}_\kappa \). Then the following results hold:
\begin{enumerate}[label=(\roman*)]
\item  If f is $\nabla$-differentiable at $t$, then $f$ is continuous at  $t.$
\item  If $f$ is continuous at $t$ and $t$ is left-scattered, then $f$ is $\nabla$-differentiable at $t$ with
\begin{equation*}
\label{rrrsder}
\nabla f(t)=\frac{f(t)-
f(\rho(t))}{\nu(t)}.
\end{equation*}
\item If $t$ is  left-dense, then $f$ is $\nabla$-differentiable at $t$ iff the limit
$$\lim_{s\to t}\frac{f(s)-f(t)}{s-t}$$
exists. Then we have
$$\nabla f(t)=\lim_{s\to t}\frac{f(s)-f(t)}{s-t}.$$
\item If $f$ is $\nabla$-differentiable at $t$, then
\begin{equation}
\label{derrho}
f(\rho(t))=f(t)-\nabla f(t)\nu(t).
\end{equation}
\end{enumerate}
\end{theorem}
We will now present the intermediate value theorem for a continuous function on a time scale.
\begin{theorem}\label{IVT}\cite{Peterson}
Assume \( f: \mathbb{T} \to \mathbb{R} \) is continuous,  \( a < b \) are points in \( \mathbb{T} \), and
$
f(a)f(b) < 0.
$
Then, there exists \( c \in [a, b)\cap\mathbb{T}\) such that either \( f(c) = 0 \) or \(f(c)  f(\sigma(c))< 0.\)
\end{theorem}

\begin{definition}\cite{Zadeh} \label{def7}
	A fuzzy set \( u \) in a universe of discourse \( U \) is represented by a function $u:U\rightarrow [0,1]$, where \( u(x) \) indicates the membership degree of \( x \) to the fuzzy set \( u \).
\end{definition}

\begin{definition}\cite{Negoita}\label{def8}
	Let $u:U\rightarrow \lbrack 0,1]$ be a fuzzy set. The $\alpha$-level sets of $u$
	are defined as
	\[
	u_{\alpha}=\left\{ x\in U:u(x)\geq \alpha\right\}
	\]%
	for $0<\alpha\leq 1$. The $0$-level set of $u$ 
	\[
	u_{0}=cl\left\{ x\in U:u(x)>0\right\} 
	\]%
	is called the support of the fuzzy set $u$. Here $cl$ denotes the closure
	of the set $u.$
\end{definition}

\begin{definition}\cite{Negoita}\label{def9}
	Let \( u: \mathbb{R} \rightarrow [0,1] \) be a fuzzy subset of the real numbers. Then \( u \) is said to be a fuzzy number if it fulfills the following criteria:

    \begin{enumerate}[label=(\roman*)]
		\item \( u \) is normal, i.e., there exists an \( x_{0} \in \mathbb{R} \) such that \( u(x_{0}) = 1 \).
		
		\item  \( u \) is quasi-concave, i.e., for all \( \lambda \in [0,1] \), $x, y\in\mathbb{R}$, \( u(\lambda x + (1-\lambda)y) \geq \min\{u(x), u(y)\} \).
		
		\item  \( u \) is upper semicontinuous on \( \mathbb{R} \), i.e., for any given \( \varepsilon > 0 \), there exists a \( \delta > 0 \) such that \( u(x) - u(x_{0}) < \varepsilon \) whenever \( |x - x_{0}| < \delta \).
		
	\item  \( u \) is compactly supported, i.e., the set \( \text{cl}\{ x \in \mathbb{R} : u(x) > 0 \} \) is compact.
    \end{enumerate}
\end{definition}
Henceforth, we will use  $\mathbb{R}_{\cal{F}}$ to denote the the space of fuzzy numbers.

\begin{definition}\label{def10**}\cite{Bedebook}
	Let $a\leq b\leq c$ be real numbers. A fuzzy number $u:%
	\mathbb{R}
	\rightarrow \lbrack 0,1]$ is called a triangular fuzzy number if%
	\[
	u(x)=\left\{ 
	\begin{array}{cc}
		\frac{x-a}{b-a}; & a\leq x\leq b, \\ 
		\frac{c-x}{c-b}; & b\leq x\leq c, \\ 
		0; & otherwise.%
	\end{array}%
	\right. 
	\]%
	Such a fuzzy number is denoted by  $u=(a,b,c)$.
\end{definition}

If $\ u=(a,b,c)$ is a triangular fuzzy number, then its $\alpha$-level sets are given by
\[
[u]_{\alpha}=[a+\alpha(b-a),c+\alpha(b-c)].
\]

\begin{theorem}\label{Bede1}\cite{Bede}
Let us consider the functions $u^{-}, u^{+}: [0,1]\to\mathbb{R}$, that satisfy the following conditions:

(i) $u^{-}(\alpha)=u^{-}_{\alpha}\in \mathbb{R}$
 is a bounded, non-decreasing, left-continuous function
in $(0,1]$ and it is right-continuous at $0.$

(ii) $u^{+}(\alpha)=u^{+}_{\alpha}\in \mathbb{R}$ is a bounded, non-increasing, left-continuous function
in $(0,1]$ and it is right-continuous at $0.$

(iii) $u^{-}_{1}\leq u^{+}_{1}.$

Then there is a fuzzy number $u\in \mathbb{R}_{\cal{F}}$ that has $u^{-}_{\alpha}, u^{+}_{\alpha}$ as endpoints of its
$\alpha$-level sets, $[u]_{\alpha}$.

\end{theorem}

\begin{definition}\cite{Stefanini}\label{def11}
	Let $u,v\in $ $\mathbb{R}_{\mathcal{F}}$. The generalized Hukuhara difference ($gH$-difference) is the fuzzy number $w$, if it exists, such that%
	$$
	u\ominus _{gH}v=w\Longleftrightarrow u=v+w \text{   or   } v=u+(-1)w. $$
 \end{definition}   

\begin{remark}\label{remark1}\cite{Stefanini}
    The conditions for the existence of \( w = u \ominus_{g_H} v \) are as follows:

\textbf{Case (i):}  
    \begin{itemize}
        \item \( w_{\alpha}^- = u_{\alpha}^- - v_{\alpha}^- \) and \( w_{\alpha}^+ = u_{\alpha}^+ - v_{\alpha}^+, \forall\alpha\in[0,1] \)
        \item \(w_{\alpha}^-\) increasing and \( w_{\alpha}^+\)  decreasing with \( w_{\alpha}^- \leq w_{\alpha}^+, \forall\alpha \in [0, 1]\)
    \end{itemize}
\textbf{Case (ii):}  
    \begin{itemize}
        \item \( w_{\alpha}^- = u_{\alpha}^+ - v_{\alpha}^+ \) and \( w_{\alpha}^+ = u_{\alpha}^- - v_{\alpha}^-, \forall\alpha\in[0,1] \)
        \item \( w_{\alpha}^- \) increasing and  \( w_{\alpha}^+ \) decreasing with \( w_{\alpha}^- \leq w_{\alpha}^+, \forall {\alpha}\in[0, 1]\).
    \end{itemize}
\end{remark}

\begin{theorem}\cite{Bede,Stefanini}\label{prop2} Let $u,v\in $ 
	$\mathbb{R}_{\mathcal{F}}
	.$ If  $u\ominus _{gH}v$ exists, then 
	\[
	[u\ominus _{gH}v] _{\alpha}=[\min
	\{u_{\alpha}^{-}-v_{\alpha}^{-}, u_{\alpha}^{+}-v_{\alpha}^{+}\}, \max
	\{u_{\alpha}^{-}-v_{\alpha}^{-}, u_{\alpha}^{+}-v_{\alpha}^{+}\}].
	\]
\end{theorem}

The Hausdorff distance for fuzzy numbers, defined below, is widely applied in fields such as fuzzy clustering, pattern recognition, and decision-making, where precise comparisons between fuzzy sets are crucial. By adapting the Hausdorff distance to account for the inherent uncertainty in fuzzy numbers, it becomes a powerful tool for evaluating the similarity or difference between fuzzy quantities in practical applications.
\begin{definition}\cite{Diamond2} \label{def12}
	The metric $D:\mathbb{R}_{\mathcal{F}}\times \mathbb{R}_{\mathcal{F}}
    \rightarrow 
	[0,\infty)$ defined by
    \begin{equation*}
    \begin{split}
	D(u,v)&=\sup_{\alpha\in \lbrack 0,1]}D_H([u]_{\alpha},[v]_{\alpha})\\
    &=\sup_{\alpha\in \lbrack 0,1]}\max \left\{ \left\vert
	u_{\alpha}^{-}-v_{\alpha}^{-}\right\vert ,\left\vert u_{\alpha}^{+}-v_{\alpha}^{+}\right\vert
	\right\}, 
    \end{split}
    \end{equation*}
	where $[u]_{\alpha}=[u_{\alpha}^{-},u_{\alpha}^{+}]$, $[v]_{\alpha}=[v_{\alpha}^{-},v_{\alpha}^{+}]$, is called the Hausdorff metric for fuzzy numbers.
\end{definition}


\begin{theorem}\cite{Diamond2} \label{prop3} Let $u, v, w, e\in \mathbb{R}_{\cal{F}}$ and $k\in\mathbb{R}$. The Hausdorff metric satisfies the followings:
	\begin{enumerate}[label=(\roman*)]
		\item $D\left( u+w,v+w\right) =D\left(u,v\right)$

  	\item $D(ku,kv) =\left\vert k\right\vert
		D\left(u,v\right)$

		\item $D\left(u+v,w+e\right) \leq D\left(
		u,w\right) +D\left(v,e\right). $
	\end{enumerate}
	
\end{theorem}

\section{Main Results}
\begin{definition}\label{defder}\cite{you}
Let $f:\mathbb{T}\to\mathbb{R}_{\cal{F}}$ be a fuzzy function, and let $t\in\mathbb{T}_{\kappa}.$ We define  the generalized Hukuhara $\nabla$-derivative of $f$ at \(t\), denoted by $\nabla_{gH} f(t),$ to be the element of $\mathbb{R}_{\cal{F}}$ (if it exists) with the property that  given any \( \varepsilon > 0\), there exists a neighborhood \( U_{\mathbb{T}}= (t- \delta, t+\delta) \cap \mathbb{T} \) of  $t$  for some \( \delta > 0 \) such that the following inequalities hold:
\begin{equation*}
D( f(t+h)\ominus
_{gH}f(\rho(t)),\nabla_{gH} f(t)(h+\nu(t)))  \leq \varepsilon |h+\nu(t)|,
\end{equation*}
\begin{equation*}
D( f(\rho(t))\ominus
_{gH}f(t-h),\nabla_{gH} f(t)(h-\nu(t))) \leq \varepsilon |h-\nu(t)|
\end{equation*}
	for all \( t-h, t+h \in U_{\mathbb{T}} \) with $0\leq h<\delta.$ 
 
     We say that $f$ is  \(\nabla_{gH} \)-differentiable at $t$ if its \(\nabla_{gH}\)-derivative exists at \(t\). Moreover, we say that $f$ is  \(\nabla_{gH} \)-differentiable on $\mathbb{T}_{\kappa}$ if its \(\nabla_{gH}\)-derivative exists at  each \(t\in\mathbb{T}_{\kappa}\).  The fuzzy function $\nabla_{gH} f:\mathbb{T}_{\kappa}\to\mathbb{R}_{\cal{F}}$ is then called the \(\nabla_{gH}\)-derivative of $f$ on $\mathbb{T}_{\kappa}.$
   \end{definition}

The definition of the  generalized Hukuhara $\nabla$-derivative can alternatively be expressed through the concept of limit as follows.
\begin{definition}
A fuzzy function $f:\mathbb{T}\to\mathbb{R}_{\cal{F}}$ is  said to be \(\nabla_{gH} \)-differentiable at
$t\in\mathbb{T}_{\kappa}$ if there exists an element  $\nabla_{gH} f(t)\in \mathbb{R}_{\cal{F}}$ such that
$$\nabla_{gH} f(t)=\lim_{h\to 0}\frac{1}{h+\nu(t)}[f(t+h)\ominus
_{gH}f(\rho(t))].$$
\end{definition}

The following theorem was presented in \cite{you} without a proof; however, we will provide the proof here.
   \begin{theorem}\cite{you}
   If the \(\nabla_{gH}\)-derivative of $f$ at $t\in\mathbb{T}_{\kappa}$  exists, then it is unique. Hence, the \(\nabla_{gH}\)-derivative is well defined.
   \end{theorem}
   
   \begin{proof}
   Let $^{1}\nabla_{gH} f(t)$ and $^{2}\nabla_{gH} f(t)$  both be the  $\nabla_{gH}$-derivative of $f$ at $t.$ Then
\begin{equation*} 
\begin{split}
   D&(^{1}\nabla_{gH} f(t),^{2}\nabla_{gH} f(t))\\&=\frac{1}{|h-\nu(t)|}D(^{1}\nabla_{gH} f(t)(h-\nu(t)),^{2}\nabla_{gH} f(t)(h-\nu(t)))\\
   &\leq \frac{1}{|h-\nu(t)|}D(^{1}\nabla_{gH} f(t)(h-\nu(t)),f(\rho(t))\ominus
_{gH}f(t-h))\\
&\,\,\,+ \frac{1}{|h-\nu(t)|} D(f(\rho(t))\ominus
_{gH}f(t-h),^{2}\nabla_{gH} f(t)(h-\nu(t)))\\
&\leq \varepsilon+\varepsilon=2\varepsilon
  \end{split}
\end{equation*}
 for all $|h-\nu(t)|\neq 0.$  Since $\varepsilon>0$ is arbitrary, we have $ D(^{1}\nabla_{gH} f(t),^{2}\nabla_{gH}f(t))=0$ and hence $^{1}\nabla_{gH} f(t)=^{2}\nabla_{gH} f(t).$
 \end{proof}
\begin{theorem}\label{derivative}\cite{you}
Assume $f:\mathbb{T}\to\mathbb{R}_{\cal{F}}$ is a  function, and let $t\in\mathbb{T}_{\kappa}$. Then, we have the followings:
\begin{enumerate}[label=(\roman*)]

\item  If $f$ is continuous at $t$ and $t$ is left-scattered, then $f$ is $\nabla_{gH}$-differentiable at $t$ with
\begin{equation}
\label{rsder}
\nabla_{gH} f(t)=\frac{f(t)\ominus
_{gH}f(\rho(t))}{\nu(t)}.
\end{equation}
\item If $t$ is  left-dense, then $f$ is $\nabla_{gH}$-differentiable at $t$ if and only if  the following two limits exist:
$$\lim_{h\to 0^+}\frac{f(t+h)\ominus
_{gH}f(t)}{h}\quad\text{and} \quad \lim_{h\to 0^+}\frac{f(t)\ominus
_{gH}f(t-h)}{h},$$
and the equality
\begin{equation}
\label{rdderlim}
\lim_{h\to 0^+}\frac{f(t+h)\ominus
_{gH}f(t)}{h}= \lim_{h\to 0^+}\frac{f(t)\ominus
_{gH}f(t-h)}{h}=\nabla_{gH} f(t)
\end{equation}
holds.
\end{enumerate}
\end{theorem}

Let $f:\mathbb{T}\to\mathbb{R}_{\cal{F}}$ be a fuzzy function. For each $\alpha\in[0,1],$ associated to $f,$
we define the family of interval valued functions $f_{\alpha}:\mathbb{T}\to{\cal K}_C$ given by
$f_{\alpha}(t)=[f(t)]_{\alpha},$ where ${\cal K}_C$ is  the family of all bounded closed intervals in $\mathbb{R}.$ For any $\alpha\in[0,1],$  we denote
$f_{\alpha}(t)=[f_{\alpha }^{-}(t), f_{\alpha }^{+}(t)],$ $t\in\mathbb{T}.$ Here, for each $\alpha\in[0,1],$  the endpoint functions $f_{\alpha }^{-},f_{\alpha }^{+}:\mathbb{T}\to\mathbb{R}$ are called the lower and upper  functions of $f,$ respectively.

In the following theorem, we focus on examining the relationship between the $\nabla_{gH}$-differentiability of a fuzzy function $f$ and the $\nabla_{gH}$-differentiability of the family of interval valued functions $f_{\alpha}.$

\begin{theorem}\label{Thmlevel}
If $f:\mathbb{T}\to\mathbb{R}_{\cal{F}}$  is $\nabla_{gH}$-differentiable at $t\in\mathbb{T}_{\kappa},$  then the family of interval valued functions $f_{\alpha}:\mathbb{T}\to{\cal K}_C$ is 
$\nabla_{gH}$-differentiable at $t$ uniformly in $\alpha\in[0,1]$. Furthermore, the $\nabla_{gH}$-derivative of $f_{\alpha}$ at $t$ is given by
$$\nabla_{gH} f_{\alpha}(t)=[\nabla_{gH} f(t)]_{\alpha }$$
for all $\alpha\in[0,1].$
\end{theorem}

\begin{proof}
Assume that $f$ is $\nabla_{gH}$-differentiable at $t.$ Let $\varepsilon>0$ be given.  Since $f$ is $\nabla_{gH}$-differentiable at $t$, there is a neighborhood $U_{\mathbb{T}}$ of $t$ such that
\begin{equation*}
\begin{split}
&D\left(\frac{1}{h+\nu(t)}[f(t+h)\ominus
_{gH}f(\rho(t))],\nabla_{gH} f(t)\right)\\
&=\underset{\alpha\in[0,1]}{\sup} D_H\left(\frac{1}{h+\nu(t)}\{[f(t+h)]_{\alpha}\ominus
_{gH}[f(\rho(t))]_{\alpha}\},[\nabla_{gH} f(t)]_{\alpha}\right)\leq\varepsilon\\
\end{split}
\end{equation*}
for all $t+h\in U_{\mathbb{T}}$ with $|h|>0.$ In view of the definition of the metric $D,$ this implies that $f_{\alpha}$ is $\nabla_{gH}$-differentiable at $t$  uniformly in $\alpha\in[0,1]$, and $$\nabla_{gH} f_{\alpha}(t)=[\nabla_{gH} f(t)]_{\alpha }$$ for all $\alpha\in[0,1].$
\end{proof}

\begin{theorem}\label{derivativeleftscattered}
Let $f:\mathbb{T}\to\mathbb{R}_{\cal{F}}$ be a function that is $\nabla_{gH}$-differentiable at a point $t\in \mathbb{T}_{\kappa}$, where $t$ is both left-scattered and right-dense. Additionally, assume that for every   \( \delta > 0 \), there exist \( h, k\in\mathbb{R}\) such that \( 0 < h, k < \delta \) and \( t + h, t+k \in \mathbb{T} \), and that the Hukuhara differences 
$$f(t+h)\ominus_H f(\rho(t))\quad\text{and}\quad
f(\rho(t))\ominus_H f(t+k)$$  are both defined. Then, $\nabla_{gH}f(t)$ is a crisp number.
\end{theorem}
\begin{proof}
Suppose the function $f$ is $\nabla_{gH}$-differentiable at a point $t\in\mathbb{T}_{\kappa},$ where $t$ is both left-scattered and right-dense.  Assuming the relevant Hukuhara differences are well-defined, then for each  \( n \in \mathbb{N} \), there exist \( h_n, h_n^{\prime} \in \mathbb{R} \) satisfying $0 < h_n<1/n,$ $ 0<h_n^{\prime}< 1/n, 
$ and such that $t+h_n, t+h_n^{\prime}\in\mathbb{T}.$  For these choices of $h_n$ and  $h_n^{\prime}$, the following inequalities hold for all \( \alpha \in [0,1] \):
\begin{equation}\label{EQ113}
\text{len}([f(t+h_n)]_{\alpha})\geq \text{len}([f(\rho(t))]_{\alpha}),
\end{equation}
\begin{equation}\label{EQ213}
\text{len}([f(t+h_h^{\prime})]_{\alpha})\leq \text{len}([f(\rho(t))]_{\alpha}).
\end{equation}

Based on this, define two sequences $\{w_n\}$ and  $\{w_n^{\prime}\}$  as:
\begin{equation*}
\begin{split}
 w_n&=\frac{1}{h_n+\nu(t)}[f(t+h_n)\ominus_{gH}f(\rho(t))],\\
 w_n^{\prime}&=\frac{1}{h_n^{\prime}+\nu(t)}[f(t+h_n^{\prime})\ominus_{gH}f(\rho(t))].
 \end{split}
 \end{equation*}

Using the inequalities   \eqref{EQ113} and \eqref{EQ213}, we can express the  $\alpha $-level sets of $w_n$ and $w_n^{\prime}$ for all \( \alpha \in [0,1] \) as:
\begin{equation}\label{EQ31}
[w_n]_{\alpha}=\frac{1}{h_n+\nu(t)}[f_{\alpha }^{-}(t+h_n)-f_{\alpha }^{-}(\rho(t)),f_{\alpha }^{+}(t+h_n)-f_{\alpha }^{+}(\rho(t))],
\end{equation}
\begin{equation}\label{EQ41}
[w_n^{\prime}]_{\alpha}=\frac{1}{h_n^{\prime}+\nu(t)}[f_{\alpha }^{+}(t+h_n^{\prime})-f_{\alpha }^{+}(\rho(t)),f_{\alpha }^{-}(t+h_n^{\prime})-f_{\alpha }^{-}(\rho(t))].
\end{equation}
By combining  the expressions \eqref{EQ31} and \eqref{EQ41} with Theorem \ref{Thmlevel},  it follows that the right-sided nabla derivatives 
 \( \nabla_+f_{\alpha}^{-}(t) \) and \(\nabla_+ f_{\alpha}^{+}(t) \) exist  for all $\alpha \in [0,1],$ and the convergence with respect to
$\alpha$  is uniform.  Furthermore, the derivatives are equal for all $\alpha \in [0,1],$  i.e.,  \(\nabla_{+} f_{\alpha}^{-}(t) = \nabla_{+}f_{\alpha}^{+}(t) \). 
 This implies that the generalized Hukuhara nabla derivative $\nabla_{gH} f(t)$ is a crisp number.
\end{proof}
In \cite{you}, the proof of the following theorem is  incomplete. Hence, we have reproven it using the result we have obtained in Theorem \ref{derivativeleftscattered}.
\begin{theorem}
Let $f:\mathbb{T}\to\mathbb{R}_{\cal{F}}$ be a  function, and let $t\in\mathbb{T}_{\kappa}$. 
 If $f$  is $\nabla_{gH}$-differentiable at $t,$ then $f$ is continuous at $t.$ 
\end{theorem}
\begin{proof}
Firstly, if $f$ is $\nabla_{gH}$-differentiable at 
 $t$ and $t$ is isolated, then it is clear that $f$ is continuous at $t.$
Now, assume that  $f$ is $\nabla_{gH}$-differentiable at 
 $t$ and $t$ is left-dense. Let $\varepsilon\in(0,1)$ be given. Define $\varepsilon^{*}=\varepsilon[1+||\nabla_{gH} f(t)||_F]^{-1}.$  Then,  $\varepsilon^{*}\in(0,1).$ Here, for $u\in\mathbb{R}_{\cal{F}},  ||u||_{F}=D(u,\tilde{0})$, where $\tilde{0}$ is the zero element of $\mathbb{R}_{\cal{F}}.$ By Definition \ref{defder},  there exists a neighborhood $U_{\mathbb{T}}$  of  $t$ such that 
\begin{equation*}
D( f(t+h)\ominus
_{gH}f(t),h\nabla_{gH} f(t))  \leq \varepsilon^* h
\end{equation*}
and
\begin{equation*}
D( f(t)\ominus
_{gH}f(t-h),h\nabla_{gH} f(t))  \leq \varepsilon^* h 
\end{equation*}
for all \( t-h, t+h \in U_{\mathbb{T}} \) with $h\geq 0.$ Hence, for all  \( t-h, t+h \in U_{\mathbb{T}}\cap(t-\varepsilon^*,t+\varepsilon^*)\) with $0\leq h<\varepsilon^*$, we have 
\begin{equation*}
\begin{split}
D(f(t-h),f(t))
&=D(f(t-h)\ominus
_{gH}f(t),\tilde{0})\\
&\leq D(f(t-h)\ominus
_{gH}f(t),(-1)h\nabla_{gH} f(t))+D((-1)h\nabla_{gH} f(t),\tilde{0})\\
&=D(f(t)\ominus
_{gH}f(t-h),h\nabla_{gH} f(t))+h D(\nabla_{gH} f(t),\tilde{0})\\
&=D(f(t)\ominus
_{gH}f(t-h),h\nabla_{gH} f(t))+h||\nabla_{gH} f(t)||_F\\
&\leq \varepsilon^* h+h||\nabla_{gH} f(t)||_F\\
&<\varepsilon^* h+\varepsilon^* ||\nabla_{gH} f(t)||_F\\
&<(1+||\nabla_{gH} f(t)||_F)\varepsilon^* \\
&=\varepsilon.
\end{split}
\end{equation*} 
Thus, we conclude that
$$\lim_{h\to 0^+}f(t-h)=f(t).$$
 If $t$ is right-scattered, then it is clear that f is right-continuous at $t.$ Otherwise, it is similarly straightforward to show that
$$\lim_{h\to 0^+}f(t+h)=f(t).$$
Therefore, $f$ is continuous at $t.$
\\Now assume that $f$ is $\nabla_{gH}$-differentiable at $t$ and $t$ is left-scattered and right-dense.  It is evident that $f$ is left continuous at $t.$ First, we consider the case where, for all $\delta>0,$ there exist $h_1,h_2$ such that $0<h_i<\delta,$ and $t+h_i\in\mathbb{T}$ for $ i=1,2,$ for which the  Hukuhara differences $f(t+h_1)\ominus_H f(\rho(t))$ and 
$f(\rho(t))\ominus_H f(t+h_2)$  exist. Then, by Theorem \ref{derivativeleftscattered}, we conclude that  $\nabla_{gH}f(t)$ is a crisp number. Let $\varepsilon\in(0,1)$ be given. Define $\varepsilon^{*}=\varepsilon[1+\nu(t)+||\nabla_{gH} f(t)||_F]^{-1}.$  Thus, $\varepsilon^{*}\in(0,1).$ By Definition \ref{defder},  there exists a neighborhood $U_{\mathbb{T}}$  of  $t$ such that  the following inequalities hold:
\begin{equation*}
D( f(t+h)\ominus
_{gH}f(\rho(t)),\nabla_{gH} f(t)(h+\nu(t)))  \leq \varepsilon^* |h+\nu(t)|,
\end{equation*}
\begin{equation*}
D( f(\rho(t))\ominus
_{gH}f(t-h),\nabla_{gH} f(t)(h-\nu(t))) \leq \varepsilon^* |h-\nu(t)|
\end{equation*}
	for all \( t-h, t+h \in U_{\mathbb{T}} \) with $h\geq 0.$ 
 Hence, we have for all  \( t-h, t+h \in U_{\mathbb{T}}\cap(t-\varepsilon^*,t+\varepsilon^*)\) with $0\leq h<\varepsilon^*:$ 
\begin{equation*}
\begin{split}
D(f(t+h),&f(\rho(t))\oplus\nu(t)\nabla_{gH} f(t))\\
&=D(f(t+h)\ominus_{gH}[f(\rho(t))\oplus\nu(t)\nabla_{gH} f(t)],\tilde{0})\\
&=D((f(t+h)\ominus_{gH}f(\rho(t)))\oplus(-1)\nu(t)\nabla_{gH} f(t),\tilde{0})\\
&=D(f(t+h)\ominus_{gH}f(\rho(t)),\nu(t)\nabla_{gH}f(t))\\
&\leq D(f(t+h)\ominus
_{gH}f(\rho(t)),\nabla_{gH} f(t)(h+\nu(t)))\\
&\,\,+D(\nabla_{gH} f(t)(h+\nu(t)),\nu(t)\nabla_{gH}f(t))\\
&=D(f(t+h)\ominus
_{gH}f(\rho(t)),\nabla_{gH} f(t)(h+\nu(t)))\\
&\,\,+D(h\nabla_{gH} f(t),0)\\
&=D(f(t+h)\ominus
_{gH}f(\rho(t)),\nabla_{gH} f(t)(h+\nu(t)))\\
&\,\,+hD(\nabla_{gH} f(t),0)\\
&=D(f(t+h)\ominus
_{gH}f(\rho(t)),\nabla_{gH} f(t)(h+\nu(t)))\\
&\,\,+h||\nabla_{gH} f(t)||_F\\
&\leq \varepsilon^*(h+\nu(t))+h||\nabla_{gH} f(t)||_F\\
&<\varepsilon^*(h+\nu(t))+\varepsilon^*||\nabla_{gH} f(t)||_F\\
&<\varepsilon^*(1+\nu(t))+\varepsilon^*||\nabla_{gH} f(t)||_F\\
&=(1+\nu(t)+||\nabla_{gH} f(t)||_F)\varepsilon^* \\
&=\varepsilon
\end{split}
\end{equation*} 
and 
\begin{equation*}
\begin{split}
D(f(t),f(\rho(t))\oplus\nu(t)\nabla_{gH} f(t))&=D(f(t)\ominus_{gH}[f(\rho(t))\oplus\nu(t)\nabla_{gH} f(t)],\tilde{0})\\
&=D((f(t)\ominus_{gH}f(\rho(t)))\oplus(-1)\nu(t)\nabla_{gH} f(t),\tilde{0})\\
&=D(f(t)\ominus_{gH}f(\rho(t)),\nu(t)\nabla_{gH} f(t))\\
&\leq\varepsilon^*\nu(t)\\
&<\varepsilon
\end{split}
\end{equation*}
 Thus, we conclude that
$$f(t)=f(\rho(t))\oplus\nu(t)\nabla_{gH} f(t),$$
which implies that
$$\lim_{h\to 0^+}f(t+h)=f(t).$$
Therefore, $f$ is continuous at $t.$
It remains to consider the cases where  there exists  $\delta>0$ such that for all 
 $0<h<\delta$ with $t+h\in\mathbb{T},$ we either have 
\begin{equation*}
\text{len}\left([f(t+h)]_{\alpha}\right)\geq\text{len}\left([f(\rho(t))]_{\alpha}\right)
\end{equation*}
or
\begin{equation*}
\text{len}\left([f(t+h)]_{\alpha}\right)\leq\text{len}\left([f(\rho(t))]_{\alpha}\right)
\end{equation*}
for all $\alpha\in[0,1].$   We will prove the first  case, and the second case follows similarly. Let $\varepsilon\in(0,1)$ be given. Define $\varepsilon^{*}=\varepsilon[1+\nu(t)+||\nabla_{gH} f(t)||_F]^{-1}.$  Then,  $\varepsilon^{*}\in(0,1).$ By Definition \ref{defder},  there exists a neighborhood $U_{\mathbb{T}}$  of  $t$ such that 
\begin{equation*}
D( f(t+h)\ominus
_{gH}f(\rho(t)),\nabla_{gH} f(t)(h+\nu(t)))  \leq \varepsilon^* |h+\nu(t)|
\end{equation*}
and
\begin{equation*}
D( f(\rho(t))\ominus
_{gH}f(t-h),\nabla_{gH} f(t)(h-\nu(t))) \leq \varepsilon^* |h-\nu(t)|
\end{equation*}
	for all \( t-h, t+h \in U_{\mathbb{T}} \) with $h\geq 0.$ 
  Thus, for all  \( t-h, t+h \in U_{\mathbb{T}}\cap(t-\varepsilon^*,t+\varepsilon^*)\) with $0<h<\varepsilon^*,$ we have 
\begin{equation*}
\begin{split}
D(f(t+h),&f(\rho(t))\oplus\nu(t)\nabla_{gH} f(t))\\
&=D(f(t+h)\ominus_{H}f(\rho(t)),\nu(t)\nabla_{gH} f(t))\\
&\leq D(f(t+h)\ominus_{H}f(\rho(t)),\nabla_{gH} f(t)(h+\nu(t)))\\&\,\,+
D(\nabla_{gH} f(t)(h+\nu(t)),
\nu(t)\nabla_{gH} f(t))\\
&=D(f(t+h)\ominus_{H}f(\rho(t)),\nabla_{gH} f(t)(h+\nu(t)))\\&\,\,+
D(h\nabla_{gH} f(t),\tilde{0}
)\\
&=D(f(t+h)\ominus_{H}f(\rho(t)),\nabla_{gH} f(t)(h+\nu(t)))\\&\,\,+
hD(\nabla_{gH} f(t),\tilde{0}
)\\
&=D(f(t+h)\ominus_{H}f(\rho(t)),\nabla_{gH} f(t)(h+\nu(t)))\\
&\,\,+h||\nabla_{gH} f(t)||_F\\
&\leq \varepsilon^*(h+\nu(t))+h||\nabla_{gH} f(t)||_F\\
&<\varepsilon^*(h+\nu(t))+\varepsilon^*||\nabla_{gH} f(t)||_F\\
&<(1+\nu(t)+||\nabla_{gH} f(t)||_F)\varepsilon^* \\
&=\varepsilon
\end{split}
\end{equation*} 
and similarly, we have
\begin{equation*}
D(f(t),f(\rho(t))\oplus\nu(t)\nabla_{gH} f(t))<\varepsilon.
\end{equation*}
Thus, we obtain
$$\lim_{h\to 0^+}f(t+h)=f(t),$$
which shows that f is continuous at $t.$
\end{proof}
\begin{theorem}\label{rhothm}
Let $f:\mathbb{T}\to\mathbb{R}_{\cal{F}}$ be a  function and let $t\in\mathbb{T}_{\kappa}$. If $f$ is $\nabla_{gH}$-differentiable at $t,$ then   the following holds:
$$f(t)=f(\rho(t))\oplus\nu(t)\nabla_{gH}f(t)$$
or
$$f(\rho(t))=f(t)\oplus(-1)\nu(t)\nabla_{gH}f(t).$$
\end{theorem}
\begin{proof}
If $\rho(t)=t,$ then $\nu(t)=0$ and it follows that
\begin{equation*}
f(t)=f(\rho(t))=f(\rho(t))\oplus\nu(t)\nabla_{gH}f(t) 
\end{equation*}
or
\begin{equation*}
f(\rho(t))=f(t)=f(t)\oplus(-1)\nu(t)\nabla_{gH}f(t).
\end{equation*}
On the other hand, if $\rho(t)<t,$ applying  Theorem \ref{derivative} (i) yields
$$f(t)\ominus_{gH}f(\rho(t))=\nu(t)\nabla_{gH}f(t).$$
From this, we conclude that
\begin{equation*}
f(t)=f(\rho(t))\oplus\nu(t)\nabla_{gH}f(t) 
\end{equation*}
or
\begin{equation*}
f(\rho(t))=f(t)\oplus(-1)\nu(t)\nabla_{gH}f(t).
\end{equation*}
This completes the proof.
\end{proof}

The following example illustrates a fuzzy number valued function that is 
$\nabla_{gH}$-differentiable, even though its lower and upper  functions 
fail to possess one-sided $\nabla$-derivatives.

\begin{example}
Let $\mathbb{T}=\left\{\frac{1}{n}: n\in\mathbb{Z}\right\} \cup \left\{\frac{\sqrt{2}}{n}: n\in\mathbb{Z}\right\} \cup \{0\}$ and let $f:\mathbb{T}\to\mathbb{R}_{\mathcal{F}}$ be defined by
$$
f(t)=\begin{cases} 
      (-2,\dfrac{t^2+t-2}{2},t^2+t), & t\in\left\{\frac{1}{n}\right\}_{n\in\mathbb{Z}} \cup \{0\},\\
      (t-2,\dfrac{t^2+t-2}{2},t^2), & t\in\left\{\frac{\sqrt{2}}{n}\right\}_{n\in\mathbb{Z}}.
   \end{cases}
$$
The corresponding lower and upper level functions are given by
$$
f_{\alpha}^-(t)=\begin{cases} 
      -2+\alpha\left(\dfrac{t^2+t+2}{2}\right), & t\in\left\{\frac{1}{n}\right\}_{n\in\mathbb{Z}} \cup \{0\},\\
     t-2+\alpha\left(\dfrac{t^2 - t + 2}{2}\right), & t\in\left\{\frac{\sqrt{2}}{n}\right\}_{n\in\mathbb{Z}},
   \end{cases}
$$
and
$$
f_{\alpha}^+(t)=\begin{cases} 
       t^2+t-\alpha\left(\dfrac{t^2+t+2}{2}\right), & t\in\left\{\frac{1}{n}\right\}_{n\in\mathbb{Z}} \cup \{0\},\\
      t^2-\alpha\left(\dfrac{t^2 - t + 2}{2}\right), & t\in\left\{\frac{\sqrt{2}}{n}\right\}_{n\in\mathbb{Z}}.
   \end{cases}
$$ 
We will now show that the right-hand nabla derivative $\nabla_+ f_{\alpha}^-(0)$ is undefined for any $\alpha \in [0,1)$. It is known that $f_{\alpha}^-$ is right $\nabla$-differentiable at $0$ if and only if the limit
\[
\nabla_+f_{\alpha}^-(0)=\lim_{h\to 0^+}\frac{f_{\alpha}^-(h)-f_{\alpha}^-(0)}{h}
\]
exists and is finite.\\
First, consider the sequence $h_n \in\left\{\frac{1}{n}\right\}_{n\in\mathbb{Z}}$.  For this choice, we have:
\[
f_{\alpha}^-(h_n) = -2 + \alpha\left(\frac{h_n^2 + h_n + 2}{2}\right), \quad\text{and}\quad f_{\alpha}^-(0) = -2 + \alpha.
\]
Thus,  the difference quotient becomes:
\[
\begin{split}
\lim_{n \to \infty} \frac{f_{\alpha}^-(h_n) - f_{\alpha}^-(0)}{h_n}
= \lim_{n \to \infty} \frac{\alpha}{h_n} \left( \frac{h_n^2 + h_n}{2} \right)
= \lim_{n \to \infty} \frac{\alpha}{2} (h_n + 1)
= \frac{\alpha}{2}.
\end{split}
\]\\
Next, take $h_n \in\left\{\frac{\sqrt{2}}{n}\right\}$. In this case:
\[
f_{\alpha}^-(h_n) = h_n - 2 + \alpha\left(\frac{h_n^2 - h_n + 2}{2}\right), \quad\text{and again}\quad  f_{\alpha}^-(0) = -2 + \alpha.
\]
Computing the corresponding difference quotient gives:
\[
\begin{split}
\lim_{n \to \infty} \frac{f_{\alpha}^-(h_n) - f_{\alpha}^-(0)}{h_n}
= \lim_{n \to \infty} \frac{h_n+ \alpha\left(\frac{h_n^2 - h_n}{2}\right)}{h_n}
= \lim_{n \to \infty} \left(1 + \frac{\alpha}{2}(h_n - 1)\right)
= 1 - \frac{\alpha}{2}.
\end{split}
\]
Since these two right-hand limits are not equal, the right-sided nabla derivative 
 $\nabla_+f_{\alpha}^-(0)$ does not exist for any $\alpha \in [0,1)$. Furthermore, a similar analysis shows that the other  nabla derivatives $\nabla_+ f_{\alpha}^+(0)$, $\nabla_- f_{\alpha}^-(0)$, and $\nabla_- f_{\alpha}^+(0)$ also fail to exist for  $\alpha \in [0,1)$. Nevertheless, the fuzzy number valued function \( f \) is $\nabla_{gH}$-differentiable at \(t=0\). To demonstrate this, we apply Theorem~\ref{derivative}~(iii). For every $\alpha \in [0,1]$,  we examine the expression:
\[
\left[\lim_{h \to 0} \frac{f(h) \ominus_{gH} f(0)}{h} \right]_\alpha.
\]
First, take the sequence $h_n \in\left\{\frac{1}{n}\right\}_{n\in\mathbb{Z}}$. Then:
\[
f(h_n) = (-2, \frac{h_n^2 + h_n - 2}{2}, h_n^2 + h_n), \quad f(0) = (-2, -1, 0).
\]
Hence, we compute
\[
\left[ \frac{f(h_n) \ominus_{gH} f(0)}{h_n} \right]_\alpha
= \left[
\min\left\{ \frac{\alpha}{2}(h_n + 1), \left(1 - \frac{\alpha}{2}\right)(h_n + 1) \right\},
\max\left\{ \frac{\alpha}{2}(h_n + 1), \left(1 - \frac{\alpha}{2}\right)(h_n + 1) \right\}
\right].
\]
Taking the limit as \( n \to \infty \), we get
\[
\lim_{n \to \infty}
\left[ \frac{f(h_n) \ominus_{gH} f(0)}{h_n} \right]_\alpha
= \left[ \frac{\alpha}{2}, 1 - \frac{\alpha}{2} \right].
\]\\
Now consider the sequence $h_n \in\left\{\frac{1}{n}\right\}_{\sqrt{2}\in\mathbb{Z}}$. In this case:
\[
f(h_n) = (h_n - 2, \frac{h_n^2 + h_n - 2}{2}, h_n^2), \quad f(0) = (-2, -1, 0).
\]
This leads to:
\[
\left[ \frac{f(h_n) \ominus_{gH} f(0)}{h_n} \right]_\alpha
= \left[
\min\left\{ 1 + \frac{\alpha}{2}(h_n - 1), h_n - \frac{\alpha}{2}(h_n - 1) \right\},
\max\left\{ 1 + \frac{\alpha}{2}(h_n - 1), h_n - \frac{\alpha}{2}(h_n - 1) \right\}
\right].
\]
Taking the limit as \( n \to \infty \), we again obtain:
\[
\lim_{n \to \infty}
\left[ \frac{f(h_n) \ominus_{gH} f(0)}{h_n} \right]_\alpha
= \left[ \frac{\alpha}{2}, 1 - \frac{\alpha}{2} \right].
\]
Combining both cases, we conclude that
\[
\left[ \lim_{h \to 0} \frac{f(h) \ominus_{gH} f(0)}{h} \right]_\alpha
= \left[ \frac{\alpha}{2}, 1 - \frac{\alpha}{2} \right],
\]
which confirms that \( f \) is indeed $\nabla_{gH}$-differentiable at \( t = 0 \) and  $\displaystyle \nabla_{gH}f(0)=\left(0,\frac{1}{2},1\right).$
\end{example}
As the preceding analysis has shown, the $\nabla_{gH}$-differentiability of a fuzzy number valued function cannot, in general, be characterized solely through the one-sided $\nabla$-differentiability of its lower and upper  functions. Such a characterization is only possible under additional assumptions that ensure the existence of these one-sided $\nabla$-derivatives.

The following theorem describes the structure of the $\alpha$-level sets of $\nabla_{gH}f(t).$ 
\begin{theorem}\label{C1}
Let $f:\mathbb{T}\to\mathbb{R}_{\cal{F}}$ be a function.  
Suppose that \( f \) is $\nabla_{gH}$-differentiable  at a point \( t \in \mathbb{T}_{\kappa} \),  
and that it is continuous on a neighborhood  
$
U_{\mathbb{T}} = (t - \delta, t + \delta) \cap \mathbb{T}
$ of $t$  
for some \( \delta > 0 \). 
 Assume further that, at every left-scattered point \(s\in U_{\mathbb{T}}\),  either the Hukuhara differences \( f(t) \ominus_H f(t + s) \) and \( f(t) \ominus_H f(t + \rho(s)) \) exist, or the differences \( f(t + s) \ominus_H f(t) \) and \( f(t+\rho(s)) \ominus_H f(t) \) exist. Then, one of the following cases holds:
\begin{enumerate}
    \item[(i)]   The functions $f_{\alpha }^{-}$ and $f_{\alpha }^{+}$ are $\nabla$-differentiable at $t$  for all $\alpha\in[0,1]$ and uniformly in $\alpha$, and either
	\begin{equation*}
		[\nabla_{gH} f(t)]_{\alpha }=\left[\nabla
		f_{\alpha }^{-}(t),\nabla
		f_{\alpha }^{+}(t)\right]
		\end{equation*}
        or
    \begin{equation*}
        [\nabla_{gH} f(t)]_{\alpha }=\left[\nabla
		f_{\alpha }^{+}(t),\nabla
		f_{\alpha }^{-}(t)\right]
    \end{equation*} 
for all $\alpha\in[0,1].$   

\item[(ii)] The left and right derivatives, $\nabla_{-} f_{\alpha }^{-}, \nabla_{+} f_{\alpha }^{-}, \nabla_{-} f_{\alpha }^{+},$ and $ \nabla_{+} f_{\alpha }^{+},$ exist at $t$ for all $\alpha\in[0,1]$ and  uniformly in $\alpha$, and satisfy  the relations
$$\nabla_{-} f_{\alpha }^{-}(t)=\nabla_{+} f_{\alpha }^{+}(t)\quad \text{and} \quad\nabla_{+} f_{\alpha }^{-}(t)=\nabla_{-} f_{\alpha }^{+}(t).$$
Furthermore, either
$$[\nabla_{gH} f(t)]_{\alpha }=[\nabla_{+} f_{\alpha }^{-}(t),\nabla_{+} f_{\alpha }^{+}(t)]=[\nabla_{-} f_{\alpha }^{+}(t),\nabla_{-} f_{\alpha }^{-}(t)]$$
or
$$[\nabla_{gH} f(t)]_{\alpha }=[\nabla_{-} f_{\alpha }^{-}(t),\nabla_{-} f_{\alpha }^{+}(t)]=[\nabla_{+} f_{\alpha }^{+}(t),\nabla_{+} f_{\alpha }^{-}(t)]$$
holds for all $\alpha\in[0,1].$ 
\end{enumerate}
\end{theorem}
\begin{proof}
Let us begin by assuming that the function $f$ is $\nabla_{gH}$-differentiable at a point $t\in\mathbb{T}_{\kappa}$ , where $t$ is a dense point in the time scale.

We now proceed by considering the following two cases:

\textbf{Case 1:} For each \( n \in \mathbb{N} \), there exist \( h_n, h_n^{\prime} \in \mathbb{R} \) such that $0 < |h_n|<1/n,$ $ 0<|h_n^{\prime}|< 1/n, 
$ and $t+h_n, t+h_n^{\prime}\in\mathbb{T},$ for which the following inequalities  are satisfied for all \( \alpha \in [0,1] \):
\begin{equation}\label{EQ1}
\text{len}\left([f(t+h_n)]_{\alpha}\right) \geq \text{len}\left([f(t)]_{\alpha}\right),
\end{equation}
\begin{equation}\label{EQ2}
\text{len}\left([f(t+h_n^{\prime})]_{\alpha}\right) \leq \text{len}\left([f(t)]_{\alpha}\right).
\end{equation}
Next, we define the sequences 
$\{w_n\}$ and $\{w_n^{\prime}\}$ as follows:
\begin{eqnarray*}\label{EQ111}
w_n &= &\frac{1}{h_n}[f(t+h_n) \ominus_{gH} f(t)],\\
w_n^{\prime} &=& \frac{1}{h_n^{\prime}}[f(t+h_n^{\prime}) \ominus_{gH} f(t)],\quad n\in \mathbb{N}.
\end{eqnarray*}
Consequently, the inequalities in \eqref{EQ1} and \eqref{EQ2} yield the following representations of the $\alpha$-level sets for all \( \alpha \in [0,1] \):
\begin{equation}\label{EQ3}
[w_n]_{\alpha} = \frac{1}{h_n}[f_{\alpha}^{-}(t+h_n) - f_{\alpha}^{-}(t), f_{\alpha}^{+}(t+h_n) - f_{\alpha}^{+}(t)],
\end{equation}
\begin{equation}\label{EQ4}
[w_n^{\prime}]_{\alpha} = \frac{1}{h_n^{\prime}}[f_{\alpha}^{+}(t+h_n^\prime) - f_{\alpha}^{+}(t), f_{\alpha}^{-}(t+h_n^\prime) - f_{\alpha}^{-}(t)].
\end{equation}
Assume that \( h_n \) and \( h_n' \) are of the same sign. Without loss of generality, we may consider the case where both
\( h_n \) and \( h_n' \) are positive. The case in which both are negative can be handled in a similar fashion.  We now aim to establish that the right-sided nabla derivative $\nabla_{+} f_{\alpha }^{-}$ exists at $t$ for every  $\alpha\in[0,1]$. To this end, assume for contradiction that $\nabla_{+} f_{\alpha }^{-}$ does not exist at $t$ for some \( \alpha_0 \in [0, 1] \).  Let \( \nabla_{gH} f_{\alpha_0}(t) = [A, B] \), where \( A, B \in \mathbb{R} \) and \( A \leq B \). Under this assumption, there exists  an \( \varepsilon > 0 \) such that, for every \( n \in \mathbb{N} \), one can find \( s_n, s_n^{\prime} \in \mathbb{R} \) satisfying
 $0 < s_n, s_n^{\prime} < 1/n,$ and $t+s_n, \, t+s_n^{\prime} \in \mathbb{T},
$
for which the following inequalities are satisfied:
\begin{equation} \label{dne}
\left| \frac{f_{\alpha_0}^{-}(t + s_n) - f_{\alpha_0}^{-}(t)}{s_n} - A \right| > \varepsilon,
\end{equation}
\begin{equation} \label{dne1}
\left| \frac{f_{\alpha_0}^{-}(t+ s_n^{\prime}) - f_{\alpha_0}^{-}(t)}{s_n^{\prime}} - B \right| > \varepsilon.
\end{equation}

We proceed by examining the situation through two separate subcases.

\noindent {\bf{Case (i):}} \( A < B \)

Given that the function $f$  is $\nabla_{gH}$-differentiable  at $t$, we have the following limits: 

\[
\min \left\{ \frac{f_{\alpha_0}^{-}(t+s_n) - f_{\alpha_0}^{-}(t)}{s_n }, \frac{f_{\alpha_0}^{+}(t+s_n) - f_{\alpha_0}^{+}(t)}{s_n }\right\} \to A
\]
and

\[
\max \left\{ \frac{f_{\alpha_0}^{-}(t+s_n^{\prime}) - f_{\alpha_0}^{-}(t)}{s_n^{\prime} }, \frac{f_{\alpha_0}^{+}(t+s_n^{\prime}) - f_{\alpha_0}^{+}(t)}{s_n^{\prime} }\right\} \to B
\]
as \(n\to\infty\).
By applying \eqref{dne} and \eqref{dne1}, 
we conclude that the following limits hold:
\[
\frac{f_{\alpha_0}^{+}(t+s_n) - f_{\alpha_0}^{+}(t)}{s_n }\to A
\]
and 
\[\frac{f_{\alpha_0}^{+}(t+s_n^{\prime}) - f_{\alpha_0}^{+}(t)}{s_n^{\prime} }\to B\]
as \( n \to \infty \). This implies that, for sufficiently large $n$, we have 
\[
\frac{f_{\alpha_0}^{+}(t+s_n) - f_{\alpha_0}^{+}(t)}{s_n }<\frac{A+B}{2}<\frac{f_{\alpha_0}^{+}(t+s_n^{\prime}) - f_{\alpha_0}^{+}(t)}{s_n^{\prime} }.
\]

We now introduce the function  $\varphi_{f_{\alpha_0}^+}:\tilde{\mathbb{T}}\to\mathbb{R}$ defined by 
$$\varphi_{f_{\alpha_0}^+}(h)=\frac{f_{\alpha_0}^{+}(t+h) - f_{\alpha_0}^{+}(t)}{h}-\frac{A+B}{2},$$  where the time scale \( \tilde{\mathbb{T}} \) is given by $\tilde{\mathbb{T}}=\{h\in\mathbb{R}: t+h\in\mathbb{T}\}.$  Since the function \( f \) is continuous on  $U_{\mathbb{T}}$, and by the definition of the Hausdorff metric, the family of functions  \(\{f_{\alpha}^{+}\}_{\alpha \in [0,1]}\) (and similarly \( \{f_{\alpha}^{-}\}_{\alpha \in [0,1]}\)) is equicontinuous on \( U_{\mathbb{T}} \).   As a result, the function \( \varphi_{f_{\alpha_0}^+} \) is continuous on a deleted neighborhood of \( 0 \), namely, on the set  \( U_{\tilde{\mathbb{T}}} = ((-\delta, \delta) \setminus \{0\}) \cap \tilde{\mathbb{T}} \).

Consider the interval \(I_n = [\min\{s_n, {s}_n^{\prime}\}, \max\{s_n,{s}_n^{\prime}\}]\cap \tilde{\mathbb{T}}\). By applying Theorem \ref{IVT}, there exists a point \(c_n \in(\min\{s_n, {s}_n^{\prime}\}, \max\{s_n,{s}_n^{\prime}\}]\cap \tilde{\mathbb{T}}\) such that either
\[
\varphi_{f_{\alpha_0}^+}(c_n) =0
\]
or
\[
\varphi_{f_{\alpha_0}^+}(c_n)\varphi_{f_{\alpha_0}^+}(\rho(c_n))<0.
\]
Observe that, under the assumptions ensuring the existence of Hukuhara differences, the second alternative is ruled out.  Consequently, we must have
\[
\frac{f_{\alpha_0}^{+}(t+c_n) - f_{\alpha_0}^{+}(t)}{c_n} = \frac{A + B}{2},
\]
with \( \lim_{n \to \infty} c_n = 0 \). However, this leads to a contradiction with the fact that the generalized Hukuhara nabla derivative satisfies
\( \nabla_{gH} f_{\alpha_0}(t) = [A, B] \) with \( A < B \). This contradiction implies that the right-sided nabla derivative $\nabla_{+} f_{\alpha }^{-}(t)$ must exist for all 
  $\alpha\in[0,1]$. By applying similar arguments, we can likewise establish the existence of the remaining one-sided derivatives $\nabla_{+} f_{\alpha }^{+}(t), \nabla_{-} f_{\alpha }^{-}(t),$ and $\nabla_{-} f_{\alpha }^{+}(t).$ 

Furthermore, from \eqref{EQ3} and \eqref{EQ4}, together with Theorem~\ref{Thmlevel}, it follows that both \( f_{\alpha}^{-} \) and \( f_{\alpha}^{+} \) are \( \nabla\)-differentiable at \( t \) for all $\alpha \in [0,1],$  with the differentiability being uniform in
 $\alpha$.  In addition, the derivatives coincide for all 
$\alpha$, that is, \( \nabla f_{\alpha}^{-}(t) = \nabla f_{\alpha}^{+}(t) \) for all \( \alpha \in [0,1] \). Consequently, the generalized Hukuhara nabla derivative \( \nabla_{gH} f(t) \) is a crisp number.

\noindent {\bf{Case (ii):}} \( A = B \)

Given that $f$ is $\nabla_{gH}$-differentiable at $t,$ we have the following limits:
\[
\min \left\{ \frac{f_{\alpha}^{-}(t+h) - f_{\alpha}^{-}(t)}{h}, \frac{f_{\alpha}^{+}(t+h) - f_{\alpha}^{+}(t)}{h } \right\} \to A
\]
and
\[
\max \left\{ \frac{f_{\alpha}^{-}(t+h) - f_{\alpha}^{-}(t)}{h}, \frac{f_{\alpha}^{+}(t+h) - f_{\alpha}^{+}(t)}{h } \right\} \to A
\]
as \( h \to 0 \). Applying the Squeeze Theorem, it follows that 
\[
\frac{f_{\alpha}^{-}(t+h) - f_{\alpha}^{-}(t)}{h} \to A
\]
and
\[
\frac{f_{\alpha}^{+}(t+h) - f_{\alpha}^{+}(t)}{h } \to A
\]
as \( h \to 0 \). Hence, the functions \( f_{\alpha}^{-} \) and \( f_{\alpha}^{+} \) are both \( \nabla \)-differentiable at \( t \) for every $\alpha \in [0,1]$, with this differentiability being uniform in $\alpha$. Moreover, both derivatives are equal to $A,$  , i.e., \( \nabla f_{\alpha}^{-}(t) = \nabla f_{\alpha}^{+}(t)=A \) for all \( \alpha \in [0,1] \). Therefore, \( \nabla_{gH} f(t) \) is a crisp number.

Suppose now that \( h_n \) and \( h_n' \) always have opposite signs for every $n.$ In other words, for each $n,$ we can consistently select either  \( h_n > 0 \) with \( h_n' < 0 \), or \( h_n < 0 \) with \( h_n' > 0 \). Under this assumption, and by invoking equations \eqref{EQ3} and \eqref{EQ4} together with Theorem \ref{Thmlevel}, we conclude that both functions \( f_{\alpha}^{-} \) and \( f_{\alpha}^{+} \) are \( \nabla\)-differentiable at \( t \) for all $\alpha \in [0,1]$,  and this differentiability is uniform in $\alpha$.  Consequently, the generalized Hukuhara 
$\nabla$-derivative of 
$f_{\alpha}$ at $t$ must take one of the following forms:
$$\nabla_{gH}f_{\alpha}(t)=[\nabla f_{\alpha }^{-}(t),\nabla f_{\alpha }^{+}(t)]$$ or $$\nabla_{gH}f_{\alpha}(t)=[\nabla f_{\alpha }^{+}(t),\nabla f_{\alpha }^{-}(t)]$$ for all $\alpha\in[0,1].$

\textbf{Case 2:} There exists a $\delta_0>0$ such that for all $h$ with
 $0<|h|<\delta_0$ and $t+h\in\mathbb{T},$ one of the following inequalities holds for all  \( \alpha \in [0,1] \):
\begin{equation}\label{EQ9}
\text{len}\left([f(t+h)]_{\alpha}\right)\geq \text{len}\left([f(t)]_{\alpha}\right)
\end{equation}
or
\begin{equation}\label{EQ10}
\text{len}\left([f(t+h)]_{\alpha}\right)\leq \text{len}\left([f(t)]_{\alpha}\right).
\end{equation}
From inequality~\eqref{EQ9}, together with Theorem~\ref{Thmlevel}, it immediately follows that the derivatives $\nabla_{-} f_{\alpha }^{-},$ $ \nabla_{+} f_{\alpha }^{-},$ $ \nabla_{-} f_{\alpha }^{+}, $ and $\nabla_{+} f_{\alpha }^{+}$ all exist at the point $t,$ with this existence being uniform in 
$\alpha.$ Additionally, these derivatives satisfy the relationships
 $$\nabla_{-} f_{\alpha }^{-}(t)=\nabla_{+} f_{\alpha }^{+}(t)\quad\text{and}\quad \nabla_{+}f_{\alpha }^{-}(t)=\nabla_{-} f_{\alpha }^{+}(t)$$ as well as
$$\nabla_{gH}f_{\alpha}(t)=[\nabla_{+} f_{\alpha }^{-}(t),\nabla_{+} f_{\alpha }^{+}(t)]=[\nabla_{-} f_{\alpha }^{+}(t),\nabla_{-} f_{\alpha }^{-}(t)]$$
for all \( \alpha \in [0,1] \).

In a similar manner, if inequality~\eqref{EQ10} holds, then by applying Theorem~\ref{Thmlevel}, the same derivatives exist and satisfy 
$$\nabla_{-} f_{\alpha }^{-}(t)=\nabla_{+} f_{\alpha }^{+}(t)\quad\text{and}\quad \nabla_{+}f_{\alpha }^{-}(t)=\nabla_{-} f_{\alpha }^{+}(t)$$ 
along with
$$\nabla_{gH}f_{\alpha}(t)=[\nabla_{-} f_{\alpha }^{-}(t),\nabla_{-} f_{\alpha }^{+}(t)]=[\nabla_{+} f_{\alpha }^{+}(t),\nabla_{+} f_{\alpha }^{-}(t)]$$
for all $\alpha\in[0,1].$ 

Following this approach, we can derive the subsequent results:

 If $f$ is $\nabla_{gH}$-differentiable at $t\in\mathbb{T}_{\kappa},$ and $t$ is left-dense and right-scattered, the following situations arise:
 
\textbf{Case 1:} The functions $f_{\alpha }^{-}$ and $f_{\alpha }^{+}$ are $\nabla$-differentiable at $t,$ uniformly with respect to $\alpha$, and they satisfy $\nabla f_{\alpha }^{-}(t)=\nabla f_{\alpha }^{+}(t)$ for all $\alpha\in[0,1].$ As a result,  $\nabla_{gH} f(t)$ is a crisp number.

\textbf{Case 2:} The functions $f_{\alpha }^{-}$ and $f_{\alpha }^{+}$ are $\nabla$-differentiable at $t$ for all $\alpha \in [0,1]$, and this differentiability is uniform with respect to $\alpha$.  Moreover, for each $\alpha\in[0,1],$  the generalized Hukuhara $\nabla$-derivative of $f_{\alpha}$ at $t$  is expressed as either
$$\nabla_{gH}f_{\alpha}(t)=[\nabla f_{\alpha }^{-}(t),\nabla f_{\alpha }^{+}(t)]$$ or $$\nabla_{gH}f_{\alpha}(t)=[\nabla f_{\alpha }^{+}(t),\nabla f_{\alpha }^{-}(t)].$$

 Assume that $f$ is $\nabla_{gH}$-differentiable at $t\in\mathbb{T}_{\kappa},$ where $t$ is right-dense and left-scattered. Given that $f$ is continuous at $t,$ and considering the definition of the Hausdorff metric, it follows that the families \( \{f_{\alpha}^{-}\}_{\alpha\in[0,1]}\) and \( \{f_{\alpha}^{+}\}_{\alpha\in[0,1]}\)  are equicontinuous at the point \( t \). According to Theorem \ref{nabladerrule} (ii), this implies that both
 \( f_{\alpha}^{-} \) and \( f_{\alpha}^{+} \) are \( \nabla\)-differentiable at \( t \),  with the differentiability holding uniformly in $\alpha.$  As a result, the generalized Hukuhara $\nabla$-derivative of $f_\alpha$ at $t$ takes one of the following forms:
 $$\nabla_{gH}f_{\alpha}(t)=[\nabla f_{\alpha }^{-}(t),\nabla f_{\alpha }^{+}(t)]$$ or $$\nabla_{gH}f_{\alpha}(t)=[\nabla f_{\alpha }^{+}(t),\nabla f_{\alpha }^{-}(t)]$$ for all $\alpha\in[0,1].$

If $f$ is $\nabla_{gH}$-differentiable at $t\in\mathbb{T}_{\kappa}$ and $t$ is isolated, then, by Theorem \ref{derivative} (i), we have $$\displaystyle
\nabla_{gH} f(t)=\frac{f(t)\ominus
_{gH}f(\rho(t))}{\nu(t)}.$$ Using the definition of the gH-difference,  for any $\alpha\in[0,1],$ we have
\begin{equation*}
[f(t)\ominus
_{gH}f(\rho(t))]_{\alpha}=[f_{\alpha }^{-}(t)-f_{\alpha }^{-}(\rho(t)),f_{\alpha }^{+}(t)-f_{\alpha }^{+}(\rho(t))]
\end{equation*}
or
\begin{equation*}
[f(t)\ominus
_{gH}f(\rho(t))]_{\alpha}=[f_{\alpha }^{+}(t)-f_{\alpha }^{+}(\rho(t)),f_{\alpha }^{-}(t)-f_{\alpha }^{-}(\rho(t))].
\end{equation*}
Multiplying by $\frac{1}{\nu(t)},$ we get
\begin{equation*}
\begin{split}
\frac{1}{\nu(t)}[f(t)\ominus
_{gH}f(\rho(t))]_{\alpha}&=\frac{1}{\nu(t)}[f_{\alpha }^{-}(t)-f_{\alpha }^{-}(\rho(t)),f_{\alpha }^{+}(t)-f_{\alpha }^{+}(\rho(t))]\\
&=\left[\frac{f_{\alpha }^{-}(t)-f_{\alpha }^{-}(\rho(t))}{\nu(t)},\frac{f_{\alpha }^{+}(t)-f_{\alpha }^{+}(\rho(t))}{\nu(t)}\right]\\
&=\left[\nabla
		f_{\alpha }^{-}(t),\nabla
		f_{\alpha }^{+}(t)\right]
\end{split}
\end{equation*}
or
\begin{equation*}
\begin{split}
\frac{1}{\nu(t)}[f(t)\ominus
_{gH}f(\rho(t))]_{\alpha}&=\frac{1}{\nu(t)}[f_{\alpha }^{+}(t)-f_{\alpha }^{+}(\rho(t)),f_{\alpha }^{-}(t)-f_{\alpha }^{-}(\rho(t))]\\
&=\left[\frac{f_{\alpha }^{+}(t)-f_{\alpha }^{+}(\rho(t))}{\nu(t)},\frac{f_{\alpha }^{-}(t)-f_{\alpha }^{-}(\rho(t))}{\nu(t)}\right]\\
&=\left[\nabla
		f_{\alpha }^{+}(t),\nabla
		f_{\alpha }^{-}(t)\right]
\end{split}
\end{equation*}
for all $\alpha\in[0,1].$ Namely,
\begin{equation*}
		\nabla_{gH}f_{\alpha}(t)=\left[\nabla
		f_{\alpha }^{-}(t),\nabla
		f_{\alpha }^{+}(t)\right]
		\end{equation*}
        or
    \begin{equation*}
        \nabla_{gH}f_{\alpha}(t)=\left[\nabla
		f_{\alpha }^{+}(t),\nabla
		f_{\alpha }^{-}(t)\right]
    \end{equation*} 
for all $\alpha\in[0,1],$ which completes the proof.

\end{proof}

The following theorem characterizes the $\nabla_{gH}$-differentiability of a fuzzy function in terms of the $\nabla$-differentiability of its endpoint functions.

\begin{theorem}
Let $f:\mathbb{T}\to\mathbb{R}_{\cal{F}}$ be a function and let $t\in\mathbb{T}_{\kappa}$.  Suppose that $f$ is continuous on a neighborhood  
$U_{\mathbb{T}} = (t - \delta, t + \delta) \cap \mathbb{T}
$ of $t$  
for some \( \delta > 0 \). 
 Assume further that, at every left-scattered point \(s\in U_{\mathbb{T}}\),  either the Hukuhara differences \( f(t) \ominus_H f(t + s) \) and \( f(t) \ominus_H f(t + \rho(s)) \) exist, or the differences \( f(t + s) \ominus_H f(t) \) and \( f(t+\rho(s)) \ominus_H f(t) \) exist. Then, $f$ is $\nabla_{gH}$-differentiable at $t$ if and only if one of the following four cases holds:

(i) The functions $f_{\alpha }^{-}$ and $f_{\alpha }^{+}$ are $\nabla$-differentiable at $t,$  uniformly with respect to $\alpha\in[0,1]$. Moreover, $\nabla
f_{\alpha }^{-}(t)$ is monotonically increasing and $\nabla
f_{\alpha }^{+}(t)$ is monotonically decreasing  as functions of $\alpha$,  and it holds that $\nabla
f_{1 }^{-}(t)\leq \nabla
f_{1 }^{+}(t).$ In this case, we have the following:
\begin{equation*}
		\left[\nabla_{gH}f(t)\right]_{\alpha}=\left[\nabla
		f_{\alpha }^{-}(t),\nabla
		f_{\alpha }^{+}(t)\right]
		\end{equation*}
for all $\alpha\in[0,1].$  

(ii) The functions $f_{\alpha }^{-}$ and $f_{\alpha }^{+}$ are $\nabla$-differentiable at $t,$ uniformly with respect to $\alpha\in[0,1]$. Moreover, $\nabla
f_{\alpha }^{+}(t)$ is monotonically increasing and $\nabla
f_{\alpha }^{-}(t)$ is monotonically decreasing  as functions of $\alpha,$ and it holds that $\nabla
f_{1 }^{+}(t)\leq \nabla
f_{1 }^{-}(t).$  In this
 case, we have the following:
\begin{equation*}
        \left[\nabla_{gH}f(t)\right]_{\alpha}=\left[\nabla
		f_{\alpha }^{+}(t),\nabla
		f_{\alpha }^{-}(t)\right]
    \end{equation*} 
for all $\alpha\in[0,1].$ 

(iii)  The one sided derivatives $\nabla_{-} f_{\alpha }^{-}, \nabla_{+} f_{\alpha }^{-}, \nabla_{-} f_{\alpha }^{+},$ and $\nabla_{+} f_{\alpha }^{+}$ exist at $t$, uniformly with respect to $\alpha\in[0,1]$.  Moreover, $\nabla_{+} f_{\alpha }^{-}(t)=\nabla_{-} f_{\alpha }^{+}(t)$ is  monotonically  increasing and $\nabla_{-} f_{\alpha }^{-}(t)=\nabla_{+} f_{\alpha }^{+}(t)$  
is  monotonically  decreasing  as functions of $\alpha$, and it holds that $\nabla_{+} f_{1}^{-}(t)\leq \nabla_{+} f_{1}^{+}(t).$ In this case, we have the following:
$$\left[\nabla_{gH}f(t)\right]_{\alpha}=[\nabla_{+} f_{\alpha }^{-}(t),\nabla_{+} f_{\alpha }^{+}(t)]=[\nabla_{-} f_{\alpha }^{+}(t),\nabla_{-} f_{\alpha }^{-}(t)]$$
for all $\alpha\in[0,1].$

(iv)  The one sided derivatives $\nabla_{-} f_{\alpha }^{-}, \nabla_{+} f_{\alpha }^{-}, \nabla_{-} f_{\alpha }^{+}, $ and $\nabla_{+} f_{\alpha }^{+}$ exist at $t$, uniformly with respect to  $\alpha\in[0,1]$. Moreover,  
$\nabla_{-} f_{\alpha }^{-}(t)=\nabla_{+} f_{\alpha }^{+}(t)$ is  monotonically  increasing and $\nabla_{+} f_{\alpha }^{-}(t)=\nabla_{-} f_{\alpha }^{+}(t)$ is  monotonically  decreasing  as functions of $\alpha,$  and it holds that $\nabla_{-} f_{1}^{-}(t)\leq \nabla_{-} f_{1}^{+}(t).$ In this case, we have the following:
$$\left[\nabla_{gH}f(t)\right]_{\alpha}=[\nabla_{-} f_{\alpha }^{-}(t),\nabla_{-} f_{\alpha }^{+}(t)]=[\nabla_{+} f_{\alpha }^{+}(t),\nabla_{+} f_{\alpha }^{-}(t)]$$
for all $\alpha\in[0,1].$
\end{theorem}
\begin{proof}
Let $f$ be $\nabla_{gH}$-differentiable at $t$. Then, $\nabla_{gH}f(t)$ is a fuzzy number, and according to Theorem \ref{C1},  one of the cases (i), (ii), (iii), or (iv) holds. 

Now assume that (i) holds. Since $f_{\alpha }^{-}$ and $f_{\alpha }^{+}$ are $\nabla$-differentiable at $t,$ uniformly with respect to $\alpha\in[0,1]$, the  limits
\begin{equation*}
 \begin{split}   
\nabla f_{\alpha }^{-}(t)&=\lim_{h\to 0 }\frac{f_{\alpha }^{-}(t+h)-f_{\alpha }^{-}(\rho(t))}{h+\nu(t)},\\
\nabla f_{\alpha }^{+}(t)&=\lim_{h\to 0 }\frac{f_{\alpha }^{+}(t+h)-f_{\alpha }^{+}(\rho(t))}{h+\nu(t)}
\end{split}
\end{equation*}
exist uniformly for all $\alpha\in[0,1].$ Let $\{h_n\}$ be a sequence that converges to 0. Then, we have
\begin{equation*}
 \begin{split}   
\nabla f_{\alpha }^{-}(t)&=\lim_{n\to\infty}\frac{f_{\alpha }^{-}(t+h_n)-f_{\alpha }^{-}(\rho(t))}{h_n+\nu(t)},\\
\nabla f_{\alpha }^{+}(t)&=\lim_{n\to\infty }\frac{f_{\alpha }^{+}(t+h_n)-f_{\alpha }^{+}(\rho(t))}{h_n+\nu(t)},
\end{split}
\end{equation*}
i.e., $\nabla f_{\alpha }^{-}(t)$ and $\nabla f_{\alpha }^{+}(t)$ are uniform limits of sequences of left continuous functions at $\alpha\in(0,1]$, and therefore they are  also  left continuous for 
$\alpha\in(0,1]$. Right continuity at 0  follows similarly.
Hence, by Theorem \ref{Bede1}, the intervals $\left[\nabla
		f_{\alpha }^{-}(t),\nabla
		f_{\alpha }^{+}(t)\right]$ define a fuzzy number. Next, assume  that $t$ is  left-dense.
We note that the following limit 
exists uniformly
\begin{equation*}
\begin{split}
&\left[\lim_{h\to 0} \frac{1}{h}[f(t+h)\ominus
_{gH}f(t)]\right]_{\alpha}\\&=\left[\lim_{h\to 0}
\frac{f_{\alpha}^{-}(t+h)-f_{\alpha}^{-}(t)}{h}, \lim_{h\to 0}\frac{f_{\alpha}^{+}(t+h)-f_{\alpha}^{+}(t)}
{h}\right]\\
&=\left[\nabla
		f_{\alpha }^{-}(t),\nabla
		f_{\alpha }^{+}(t)\right]
\end{split}
\end{equation*}
and it is a fuzzy number. Now, assume that $t$ is left-scattered.  Again, the following limit exists uniformly
\begin{equation*}
\begin{split}
&\left[\lim_{h\to 0} \frac{1}{h+\nu(t)}[f(t+h)\ominus
_{gH}f(\rho(t))]\right]_{\alpha}\\&=\left[\lim_{h\to 0}
\frac{f_{\alpha}^{-}(t+h)-f_{\alpha}^{-}(\rho(t))}{h+\nu(t)},\lim_{h\to 0} 
\frac{f_{\alpha}^{+}(t+h)-f_{\alpha}^{+}(\rho(t))}{h+\nu(t)}
\right]\\
&=\left[\frac{f_{\alpha}^{-}(t)-f_{\alpha}^{-}(\rho(t))}{\nu(t)}, \frac{f_{\alpha}^{+}(t)-f_{\alpha}^{+}(\rho(t))}
{\nu(t)}   \right]\\
&=\left[\nabla
		f_{\alpha }^{-}(t),\nabla
		f_{\alpha }^{+}(t)\right]
\end{split}
\end{equation*} 
and it is a fuzzy number. Hence, we conclude that  $f$ is
$\nabla_{gH}$-differentiable at $t$ and
\begin{equation*}
\left[\nabla_{gH}f(t)\right]_{\alpha}=\left[\nabla
		f_{\alpha }^{-}(t),\nabla
		f_{\alpha }^{+}(t)\right]
\end{equation*}
for all $\alpha\in[0,1].$ The case (ii) is analogous.

Next, suppose that (iii) holds. It can easily be verifed as in the proof of case (i) that  the pairs of functions $\nabla_{+} f_{\alpha }^{-}(t)$ and $\nabla_{+} f_{\alpha }^{+}(t),$ 
as well as  $\nabla_{-} f_{\alpha }^{+}(t)$ and $\nabla_{-} f_{\alpha }^{-}(t)$, are left-continuous on the interval $(0,1]$ and  right-continuous at $\alpha=0.$ 
 Hence, by Theorem \ref{Bede1}, the intervals
 $[\nabla_{+} f_{\alpha }^{-}(t),\nabla_{+} f_{\alpha }^{+}(t)]$ and  $[\nabla_{-} f_{\alpha }^{+}(t),\nabla_{-} f_{\alpha }^{-}(t)]$  define  fuzzy numbers. Now, assume  that $t$ is  dense. We observe that the following limits exist uniformly with respect to $\alpha\in[0,1]:$
 \begin{equation*}
\begin{split}
&\left[\lim_{h\to 0^+} \frac{1}{h}[f(t+h)\ominus
_{gH}f(t)]\right]_{\alpha}\\&=\left[\lim_{h\to 0^+}
\frac{f_{\alpha}^{-}(t+h)-f_{\alpha}^{-}(t)}{h}, \lim_{h\to 0^+}\frac{f_{\alpha}^{+}(t+h)-f_{\alpha}^{+}(t)}
{h}\right]\\
&=\left[\nabla_+
		f_{\alpha }^{-}(t),\nabla_+
		f_{\alpha }^{+}(t)\right],
\end{split}
\end{equation*}
\begin{equation*}
\begin{split}
&\left[\lim_{h\to 0^-} \frac{1}{h}[f(t+h)\ominus
_{gH}f(t)]\right]_{\alpha}\\&=\left[\lim_{h\to 0^-}
\frac{f_{\alpha}^{+}(t+h)-f_{\alpha}^{+}(t)}{h}, \lim_{h\to 0^-}\frac{f_{\alpha}^{-}(t+h)-f_{\alpha}^{-}(t)}
{h}\right]\\
&=\left[\nabla_-
		f_{\alpha }^{+}(t),\nabla_-
		f_{\alpha }^{-}(t)\right]
\end{split}
\end{equation*}
and  each of these intervals represents a fuzzy number. Using the assumption $\nabla_+
		f_{\alpha }^{-}(t)=\nabla_-
		f_{\alpha }^{+}(t)$ and $\nabla_+
		f_{\alpha }^{+}(t)=\nabla_-
		f_{\alpha }^{-}(t)$, we
        conclude that $f$ is $\nabla_{gH}$-differentiable at $t$. Furthermore, 
$$\left[\nabla_{gH}f(t)\right]_{\alpha}=[\nabla_{+} f_{\alpha }^{-}(t),\nabla_{+} f_{\alpha }^{+}(t)]=[\nabla_{-} f_{\alpha }^{+}(t),\nabla_{-} f_{\alpha }^{-}(t)]$$
for all $\alpha\in[0,1].$ Now, assume that $t$ is left-scattered. In this case,  the functions $f_{\alpha }^{-}$ and $f_{\alpha }^{+}$ are $\nabla$-differentiable at $t,$
uniformly with respect to $\alpha\in[0,1].$ Then, the following limits also exist uniformly
\begin{equation*}
\begin{split}
&\left[\lim_{h\to 0^+} \frac{1}{h+\nu(t)}[f(t+h)\ominus
_{gH}f(\rho(t))]\right]_{\alpha}\\&=\left[\lim_{h\to 0^+}
\frac{f_{\alpha}^{-}(t+h)-f_{\alpha}^{-}(\rho(t))}{h+\nu(t)},\lim_{h\to 0^+} 
\frac{f_{\alpha}^{+}(t+h)-f_{\alpha}^{+}(\rho(t))}{h+\nu(t)}
\right]\\
&=\left[\nabla_+
		f_{\alpha }^{-}(t),\nabla_+
		f_{\alpha }^{+}(t)\right],
\end{split}
\end{equation*}  
\begin{equation*}
\begin{split}
&\left[\lim_{h\to 0^-} \frac{1}{h+\nu(t)}[f(t+h)\ominus
_{gH}f(\rho(t))]\right]_{\alpha}\\&=\left[\lim_{h\to 0^-}
\frac{f_{\alpha}^{+}(t+h)-f_{\alpha}^{+}(\rho(t))}{h+\nu(t)},\lim_{h\to 0^-} 
\frac{f_{\alpha}^{-}(t+h)-f_{\alpha}^{-}(\rho(t))}{h+\nu(t)}
\right]\\
&=\left[\nabla_-
		f_{\alpha }^{+}(t),\nabla_-
		f_{\alpha }^{-}(t)\right]
\end{split}
\end{equation*} 
and both expressions represent fuzzy numbers. Using the same assumptions
on the derivatives as above, we again deduce that $f$ is
$\nabla_{gH}$-differentiable at $t.$  Furthermore, in view of the $\nabla$-differentiability of  $f_{\alpha }^{-}$ and $f_{\alpha }^{+}$, it follows that $\nabla_{gH}f(t)$ is a crisp number. The case (iv) is analogous.
 \end{proof}

We will now derive the $\nabla_{gH}$-derivative of the sum 
of two $\nabla_{gH}$-differentiable fuzzy functions on time scales.  To begin, we first define the concepts of
 (i) and (ii)-$\nabla_{gH}$-differentiability as follows.
\begin{definition}
	Let $f:\mathbb{T}\to\mathbb{R}_{\cal{F}}$ and $t\in\mathbb{T}_{\kappa}$ with $f_{\alpha }^{-}$ and $f_{\alpha }^{+}$ both $\nabla$-differentiable at $t.$ 
	We say that 
 \begin{enumerate}[label=(\roman*)]   
	\item  $f$ is (i)-$\nabla_{gH}$-differentiable at $t$ if
    $$\left[\nabla_{gH}f(t)\right]_{\alpha}=\left[\nabla
		f_{\alpha }^{-}(t),\nabla
		f_{\alpha }^{+}(t)\right],\forall\alpha\in[0,1]$$

  \item $f$ is (ii)-$\nabla_{gH}$-differentiable at $t$ if
  
  $$\left[\nabla_{gH}f(t)\right]_{\alpha}=\left[\nabla
		f_{\alpha }^{+}(t),\nabla
		f_{\alpha }^{-}(t)\right],\forall\alpha\in[0,1].$$
\end{enumerate}
\end{definition}
\begin{remark}
As demonstrated in Theorem \ref{C1}, it is possible for $f : \mathbb{T} \to \mathbb{R}_{\mathcal{F}}$ to be $\nabla_{gH}$-differentiable at $t \in \mathbb{T}_{\kappa}$, yet not be (i)-$\nabla_{gH}$-differentiable or (ii)-$\nabla_{gH}$-differentiable.
\end{remark}
\begin{theorem}
Let  $f,g:\mathbb{T}\to\mathbb{R}_{\cal{F}}$  be $\nabla_{gH}$-differentiable  at $t\in\mathbb{T}_{\kappa}$. Then, the sum   $f\oplus g:\mathbb{T}\to\mathbb{R}_{\cal{F}}$ is $\nabla_{gH}$-differentiable at $t$  with
\begin{equation}\label{PRD}
\nabla_{gH}(f\oplus g)(t)=\nabla_{gH}f(t)\oplus\nabla_{gH}g(t)
\end{equation}
provided that $f$ and $g$ are simultaneously 
(i)-$\nabla_{gH}$-differentiable or (ii)-$\nabla_{gH}$-differentiable at $t$.
\end{theorem}
\begin{proof}
Suppose that $f$ and $g$ are simultaneously 
(i)-$\nabla_{gH}$-differentiable at $t \in \mathbb{T}_{\kappa}$. Then, for $h \geq 0$, by the definition of the $gH$-difference, we obtain $f(t+h)=f(\rho(t))\oplus(f(t+h)\ominus_{gH}f(\rho(t)))$ and $g(t+h)=g(\rho(t))\oplus(g(t+h)\ominus_{gH}g(\rho(t)))$, and thus 
$$(f\oplus g)(t+h)=
(f\oplus g)(\rho(t))\oplus(f(t+h)\ominus_{gH}f(\rho(t)))\oplus(g(t+h)\ominus_{gH}g(\rho(t))).$$
That is,
\[
(f \oplus g)(t+h) \ominus_{gH} (f \oplus g)(\rho(t)) = (f(t+h) \ominus_{gH} f(\rho(t))) \oplus (g(t+h) \ominus_{gH} g(\rho(t))).
\]
For $h > 0$, we have
$f(\rho(t))=f(t-h)\oplus(f(\rho(t))\ominus_{gH}f(t-h))$ and $g(\rho(t))=g(t-h)\oplus(g(\rho(t))\ominus_{gH}g(t-h))$, 
leading to
\[
(f \oplus g)(\rho(t)) = (f \oplus g)(t-h)\oplus(f(\rho(t)) \ominus_{gH} f(t-h)) \oplus (g(\rho(t)) \ominus_{gH} g(t-h)).
\]
Namely,
\[
(f \oplus g)(\rho(t))\ominus_{gH}(f \oplus g)(t-h)=(f(\rho(t)) \ominus_{gH} f(t-h)) \oplus (g(\rho(t)) \ominus_{gH} g(t-h)).
\]

Since $f$ and $g$ are $\nabla_{gH}$-differentiable at $t$, for any given $\varepsilon > 0$, there exist neighborhoods $U_{\mathbb{T}}$ and $V_{\mathbb{T}}$ of $t$ such that
\[
D\left(f(t+h) \ominus_{gH} f(\rho(t)), \nabla_{gH} f(t)(h + \nu(t))\right) \leq \frac{\varepsilon}{2} |h + \nu(t)|
\]
and
\[
D\left(f(\rho(t)) \ominus_{gH} f(t-h), \nabla_{gH} f(t)(h - \nu(t))\right) \leq \frac{\varepsilon}{2} |h - \nu(t)|
\]
for all \( t - h, t + h \in U_{\mathbb{T}} \) with $h \geq 0$, and
\[
D\left(g(t+h) \ominus_{gH} g(\rho(t)), \nabla_{gH} g(t)(h + \nu(t))\right) \leq \frac{\varepsilon}{2} |h + \nu(t)|
\]
and
\[
D\left(g(\rho(t)) \ominus_{gH} g(t-h), \nabla_{gH} g(t)(h - \nu(t))\right) \leq \frac{\varepsilon}{2} |h - \nu(t)|
\]
for all \( t - h, t + h \in V_{\mathbb{T}} \) with $h \geq 0$. Let $W_{\mathbb{T}} = U_{\mathbb{T}} \cap V_{\mathbb{T}}$. Then, for all \( t - h, t + h \in W_{\mathbb{T}} \) with $h \geq 0$, we have
\[
\begin{split}
& D\left( (f \oplus g)(t+h) \ominus_{gH} (f \oplus g)(\rho(t)), (\nabla_{gH} f(t) \oplus \nabla_{gH} g(t))(h + \nu(t)) \right) \\
& = D\left( (f(t+h) \ominus_{gH} f(\rho(t))) \oplus (g(t+h) \ominus_{gH} g(\rho(t))), \right. \\
& \quad \quad \left. \nabla_{gH} f(t)(h + \nu(t)) \oplus \nabla_{gH} g(t)(h + \nu(t)) \right) \\
& \leq D\left( f(t+h) \ominus_{gH} f(\rho(t)), \nabla_{gH} f(t)(h + \nu(t)) \right) \\
& \quad + D\left( g(t+h) \ominus_{gH} g(\rho(t)), \nabla_{gH} g(t)(h + \nu(t)) \right) \\
& \leq \frac{\varepsilon}{2} |h + \nu(t)| + \frac{\varepsilon}{2} |h + \nu(t)| \\
& = \varepsilon |h + \nu(t)|.
\end{split}
\]
Similarly,
\[
\begin{split}
& D\left( (f \oplus g)(\rho(t)) \ominus_{gH} (f \oplus g)(t-h), (\nabla_{gH} f(t) \oplus \nabla_{gH} g(t))(h - \nu(t)) \right) \\
& = D\left( (f(\rho(t)) \ominus_{gH} f(t-h)) \oplus (g(\rho(t)) \ominus_{gH} g(t-h)), \right. \\
& \quad \quad \left. \nabla_{gH} f(t)(h - \nu(t)) \oplus \nabla_{gH} g(t)(h - \nu(t)) \right) \\
& \leq D\left( f(\rho(t)) \ominus_{gH} f(t-h), \nabla_{gH} f(t)(h - \nu(t)) \right) \\
& \quad + D\left( g(\rho(t)) \ominus_{gH} g(t-h), \nabla_{gH} g(t)(h - \nu(t)) \right) \\
& \leq \frac{\varepsilon}{2} |h - \nu(t)| + \frac{\varepsilon}{2} |h - \nu(t)| \\
& = \varepsilon |h - \nu(t)|.
\end{split}
\]
Therefore, $f+g$ is $\nabla_{gH}$-differentiable at $t$, and we have
$
\nabla_{gH} (f \oplus g)(t) = \nabla_{gH} f(t) \oplus \nabla_{gH} g(t).$
Similarly, if $f$ and $g$ are 
(ii)-$\nabla_{gH}$-differentiable at $t$, it can be shown that $f+g$ is $\nabla_{gH}$-differentiable at $t$, and equation \eqref{PRD} holds as well.
\end{proof}

In the following theorems, we restrict the fuzzy number valued function $f$ to a real function and study the generalized Hukuhara nabla derivative of the product function $fg.$ For simplicity, we define $v(t,h)=g(t+h)\ominus_{gH} g(\rho(t))$, $U(t,h)=f(t+h)-f(\rho(t)),$ and $W(t,h)=f(t+h)g(t+h)\ominus_{gH}f(\rho(t))g(\rho(t)).$ 
\begin{theorem}\label{product1}
Let $f:\mathbb{T}\to{\mathbb{R}}$ be continuously differentiable and suppose $g:\mathbb{T}\to\mathbb{R}_{\cal{F}}$  is (i)-$\nabla_{gH}$-differentiable. If $f(t)\nabla f(t)>0,$ then 
$fg:\mathbb{T}\to\mathbb{R}_{\cal{F}}$ is $\nabla_{gH}$-differentiable  and the formula
\begin{equation}\label{productrule}
\begin{split}
\nabla_{gH}(fg)(t)&=\nabla f(t)g(\rho(t))\oplus f(t)\nabla_{gH}g(t)\\
&=f(\rho(t))\nabla_{gH}g(t)\oplus\nabla f(t)g(t)
\end{split}
\end{equation}
holds.
\end{theorem}
\begin{proof}
First, let us consider the case where $h>0.$ Since $g$ is (i)-$\nabla_{gH}$-differentiable, we have
$g(t+h)=g(\rho(t))\oplus v(t,h).$ Furthermore, since  $f(t)\nabla f(t)>0,$ it follows that $f(\rho(t))$ and $U(t,h)$ have the same sign. Therefore, we have
\begin{equation*}
\begin{split}
f(t+h)g(t+h)&=f(\rho(t))g(\rho(t))\oplus U(t,h)g(\rho(t))\oplus f(\rho(t))v(t,h)\\&\,\,\oplus U(t,h)v(t,h)\\
&=f(\rho(t))g(\rho(t))\oplus U(t,h)g(\rho(t))\oplus f(t+h)v(t,h).
\end{split}
\end{equation*}
Thus,
\begin{equation}\label{W1}
W(t,h)=U(t,h)g(\rho(t))\oplus f(t+h)v(t,h).
\end{equation}
On the other hand, if $h<0,$  by the  (i)-$\nabla_{gH}$-differentiability of $g$, we have
$g(\rho(t))=g(t+h)\oplus (-1)v(t,h).$ And, in view of $f(t)\nabla f(t)>0,$ we see that  $f(t+h)$ and $(-1)U(t,h)$ have the same sign. Then,
\begin{equation*}
\begin{split}
f(\rho(t))g(\rho(t))&=f(t+h)g(t+h)\oplus(-1)U(t,h)g(t+h)\oplus (-1)f(t+h)v(t,h)\\&\,\,\oplus U(t,h)v(t,h)\\
&=f(t+h)g(t+h)\oplus (-1)U(t,h)g(t+h)\oplus (-1)f(\rho(t))v(t,h).
\end{split}
\end{equation*}
Hence,
\begin{equation}\label{W2}
W(t,h)=U(t,h)g(t+h)\oplus f(\rho(t))v(t,h).
\end{equation}
In view of \eqref{W1}, \eqref{W2}, and the differentiability and continuity of $f$ and $g$, we get
\begin{equation}\label{PR1}
\begin{split}
\lim_{h\to 0^+}\frac{1}{h+\nu(t)}W(t,h)&=\lim_{h\to 0^+}\left(\frac{1}{h+\nu(t)}U(t,h)\right)g(\rho(t))\\&\,\,\oplus f(t) \lim_{h\to 0^+}\frac{1}{h+\nu(t)}v(t,h)\\
&=\nabla f(t)g(\rho(t))\oplus f(t)\nabla_{gH}g(t)\\
\end{split}   
\end{equation}
and
\begin{equation}
\begin{split}\label{PR2}
\lim_{h\to 0^-}\frac{1}{h+\nu(t)}W(t,h)&=\lim_{h\to 0^-}\left(\frac{1}{h+\nu(t)}U(t,h)\right)g(t)\\&\,\,\oplus f(\rho(t)) \lim_{h\to 0^-}\frac{1}{h+\nu(t)}v(t,h)\\
&=\nabla f(t)g(t)\oplus f(\rho(t))\nabla_{gH}g(t).
\end{split}   
\end{equation}
Now, using the  (i)-$\nabla_{gH}$-differentiability of $g$, the condition $f(t)\nabla f(t)>0,$ Theorem \ref{nabladerrule}  (ii), and Theorem \ref{rhothm}, it is straightforward to verify that the following equality holds:
\begin{equation}
\begin{split}\label{PR3}
\nabla f(t)g(\rho(t))\oplus f(t)\nabla_{gH}g(t)&=\nabla f(t)g(\rho(t))\oplus (f(\rho(t))+\nu(t)\nabla f(t))\nabla_{gH}g(t)\\
&=\nabla f(t)g(\rho(t))\oplus f(\rho(t))\nabla_{gH}g(t)+\nu(t)\nabla f(t)\nabla_{gH}g(t)\\
&=\nabla f(t)(g(\rho(t))\oplus \nu(t)\nabla_{gH}g(t))\oplus f(\rho(t))\nabla_{gH}g(t)\\
&=\nabla f(t)g(t)\oplus f(\rho(t))\nabla_{gH}g(t).
\end{split}   
\end{equation}
Then,  by utilizing   \eqref{PR1}, \eqref{PR2}, and \eqref{PR3}, we can conclude that the function $fg$ is $\nabla_{gH}$-differentiable at $t,$
and its derivative is as stated in \eqref{productrule}.
\end{proof}

In the special case where the function $g$ is interval valued, we also derive the following result.

\begin{theorem}\label{product11}
Let $f:\mathbb{T}\to{\mathbb{R}}$ be continuously differentiable and suppose $g:\mathbb{T}\to\cal{K}_C$  is (i)-$\nabla_{gH}$-differentiable. Suppose also that 
 $f(t)\nabla f(t)<0$  and $\text{len}(f(t)g(t))$ is increasing. Then 
$fg:\mathbb{T}\to\cal{K}_C$ is $\nabla_{gH}$-differentiable  and 
\begin{equation}\label{productrule1}
\begin{split}
\nabla_{gH}(fg)(t)\oplus (-1)\nabla f(t)g(\rho(t))=f(t)\nabla_{gH}g(t)\\
\end{split}
\end{equation}
and if $\text{len}(f(t)g(t))$ is decreasing, then $fg$ is  $\nabla_{gH}$-differentiable and 
\begin{equation}
\begin{split}\label{productrule111}
\nabla_{gH}(fg)(t)\oplus (-1) f(t)\nabla_{gH}g(t)=\nabla f(t)g(\rho(t)).\\
\end{split}
\end{equation}
\end{theorem}
\begin{proof}
First, we  consider the case where $h>0.$ By the (i)-$\nabla_{gH}$-differentiability of
$g$, we know that
$g(t+h)=g(\rho(t))\oplus v(t,h).$ Since $f(t)\nabla f(t)<0,$ we conclude that $f(t+h)$ and $(-1)U(t,h)$ have the same sign. Therefore, we can write
$$(f(t+h)+(-1)U(t,h))g(t+h)=f(\rho(t))(g(\rho(t))\oplus v(t,h)),$$
which yields in
\begin{equation}\label{E1}
f(t+h)g(t+h)\oplus(-1)U(t,h)g(t+h)=f(\rho(t))g(\rho(t))\oplus f(\rho(t)) v(t,h).
\end{equation}
Assuming that $\text{len}(f(t)g(t))$ is increasing, we have
\begin{equation}\label{E2}
f(t+h)g(t+h)=f(\rho(t))g(\rho(t))\oplus W(t,h).
\end{equation}
By  substituting \eqref{E2} into \eqref{E1}, we get
\begin{equation*}\label{E3}
W(t,h)\oplus (-1)U(t,h)g(t+h)=f(\rho(t)) v(t,h).
\end{equation*} 
Then, considering the differentiability of $f$ and $g$, along with the continuity of $g$, we obtain
\begin{equation*}\label{E4}
\begin{split}
\lim_{h\to 0^+}\frac{1}{h+\nu(t)}W(t,h)&\oplus (-1)\lim_{h\to 0^+}\left(\frac{1}{h+\nu(t)}U(t,h)\right)g(t)\\&=f(\rho(t)) \lim_{h\to 0^+}\frac{1}{h+\nu(t)}v(t,h),
\end{split}
\end{equation*} 
and thus, we conclude that
\begin{equation}
\begin{split}\label{limres}
\lim_{h\to 0^+}\frac{1}{h+\nu(t)}W(t,h)\oplus (-1)\nabla f(t)g(t)=f(\rho(t))\nabla_{gH}g(t).\\
\end{split}
\end{equation}
Next,  using the  (i)-$\nabla_{gH}$-differentiability of $g$, the condition $f(t)\nabla f(t)<0,$ Theorem \ref{nabladerrule}  (ii), and Theorem \ref{rhothm}, we observe that \eqref{limres} is equivalent to
\begin{equation}\label{E45}
\begin{split}
\lim_{h\to 0^+}\frac{1}{h+\nu(t)}W(t,h)\oplus (-1)\nabla f(t)g(\rho(t))=f(t)\nabla_{gH}g(t).\\
\end{split}
\end{equation}
We now consider the case  when $h<0.$ Since  $g$
 is (i)-$\nabla_{gH}$-differentiable, we  have that
$g(\rho(t))=g(t+h)\oplus (-1)v(t,h).$  Furthermore, since $f(t)\nabla f(t)<0,$ we also know that $f(\rho(t))$ and $U(t,h)$ have the same sign.  Therefore,
 we can write
$$f(t+h)(g(t+h)\oplus (-1)v(t,h))=(U(t,h)+f(\rho(t)))g(\rho(t)),$$
which leads to
\begin{equation}\label{E11}
f(t+h)g(t+h)\oplus(-1)f(t+h)v(t,h)=U(t,h)g(\rho(t))\oplus f(\rho(t))g(\rho(t)) .
\end{equation}
Note that our assumption  $\text{len}(f(t)g(t))$ is increasing implies 
\begin{equation}\label{E22}
f(\rho(t))g(\rho(t))=f(t+h)g(t+h)\oplus(-1)W(t,h).
\end{equation}
Now, by substituting \eqref{E22} into \eqref{E11}, we observe that
\begin{equation*}\label{E33}
(-1)f(t+h)v(t,h)=U(t,h)g(\rho(t))\oplus(-1)W(t,h),
\end{equation*} 
namely,
\begin{equation*}\label{E34}
f(t+h)v(t,h)=(-1)U(t,h)g(\rho(t))\oplus W(t,h).
\end{equation*} 
Then, again  by the differentiability of $f$ and $g$ and the continuity of $f,$ we have
\begin{equation*}
\begin{split}
f(t)\lim_{h\to 0^-}\frac{1}{h+\nu(t)}v(t,h)&=(-1)\lim_{h\to 0^-}\left(\frac{1}{h+\nu(t)}U(t,h)\right)g(\rho(t))\\
&\oplus
\lim_{h\to 0^-}\frac{1}{h+\nu(t)}W(t,h)
\end{split}
\end{equation*} 
and hence,
\begin{equation}\label{E44}
f(t)\nabla_{gH}g(t)=(-1)\nabla f(t)g(\rho(t))\oplus
\lim_{h\to 0^-}\frac{1}{h+\nu(t)}W(t,h).
\end{equation}
Now,  from \eqref{E45} and \eqref{E44}, we deduce that  $fg$ is $\nabla_{gH}$-differentiable at $t$,
and the derivative satisfies the expression given in \eqref{productrule1},  as asserted. In a similar manner, it can be shown that the equation in \eqref{productrule111}
also holds.
\end{proof}

Similarly, when g is (ii)-$\nabla_{gH}$-differentiable, we  obtain the following theorems, the proofs of which are omitted.
\begin{theorem}\label{product2}
Let $f:\mathbb{T}\to{\mathbb{R}}$ be continuously differentiable and suppose $g:\mathbb{T}\to\mathbb{R}_{\cal{F}}$  is (ii)-$\nabla_{gH}$-differentiable. If $f(t)\nabla f(t)<0,$ then 
$fg:\mathbb{T}\to\mathbb{R}_{\cal{F}}$ is $\nabla_{gH}$-differentiable  and the formula
\begin{equation*}
\begin{split}
\nabla_{gH}(fg)(t)&=\nabla f(t)g(\rho(t))\oplus f(t)\nabla_{gH}g(t)\\
&=f(\rho(t))\nabla_{gH}g(t)\oplus\nabla f(t)g(t)
\end{split}
\end{equation*}
holds.
\end{theorem}
\begin{theorem}\label{product11222}
Let $f:\mathbb{T}\to{\mathbb{R}}$ be continuously differentiable and suppose $g:\mathbb{T}\to\cal{K}_C$  is (ii)-$\nabla_{gH}$-differentiable. Suppose also that 
 $f(t)\nabla f(t)>0$  and $\text{len}(f(t)g(t))$ is increasing. Then 
$fg:\mathbb{T}\to\cal{K}_C$ is $\nabla_{gH}$-differentiable  and 
\begin{equation}\label{productrule12}
\begin{split}
\nabla_{gH}(fg)(t)\oplus (-1) f(\rho(t))\nabla_{gH}g(t)=\nabla f(t)g(t)\\
\end{split}
\end{equation}
and if $\text{len}(f(t)g(t))$ is decreasing, then $fg$ is  $\nabla_{gH}$-differentiable and 
\begin{equation}
\begin{split}\label{productrule1112}
\nabla_{gH}(fg)(t)\oplus (-1)\nabla f(t)g(t)=f(\rho(t))\nabla_{gH}g(t).\\
\end{split}
\end{equation}
\end{theorem}

\section{Conclusion}

 In this paper, we have  presented an in-depth study of the generalized Hukuhara nabla derivative for fuzzy functions defined on time scales. We have derived several characterizations of fuzzy functions that are generalized Hukuhara nabla differentiable on time scales, utilizing the nabla differentiability of their endpoint functions. These characterizations complete and extend earlier results for fuzzy number valued functions in the real domain to the broader framework of time scales. Furthermore, we have extended the product rule, initially established for interval valued functions in the real domain, to fuzzy number valued functions on time scales.

This work contributes to the generalization of fuzzy calculus on time scales, providing new insights into the behavior and differentiability of fuzzy functions. The established characterizations and the extended product rule offer a more comprehensive understanding of fuzzy calculus within this framework, which can be applied in various mathematical and engineering contexts where time scales are involved.
\section*{conflict of interest}
The authors declare that there are no conflicts of interest related to this work.


\begin{thebibliography}{1}

\bibitem{Hilger1} Hilger, S., Ein Makettenkalkuls mit Anwendung auf Zentrumsmannigfaltigkeiten. Ph.D. Thesis, Universität Würzburg, Würzburg, Germany, 1988.

\bibitem{Peterson} M. Bohner, A. Peterson, Dynamic Equations on Time Scales: An Introduction with Applications, Birkhäuser, Boston, 2001.

\bibitem{Peterson1} M. Bohner, A. Peterson, Advances in Dynamic Equations on Time Scales, Birkhäuser, Boston, 2002.

\bibitem{Guseinov} F. Merdivenci Atici, G. Sh. Guseinov, On Green’s functions and positive solutions for boundary value problems on time scales, \textit{Journal of Computational and Applied Mathematics}, \textbf{141} (2002), 75–99.

\bibitem{Anderson} D. Anderson, J. Bullock, L. Erbe, A. Peterson, H. Tran, Nabla dynamic equations, in: Advances in Dynamic Equations on Time Scales, Boston, MA: Birkhäuser Boston, 2003.

\bibitem{Zadeh} L. A. Zadeh, Fuzzy Sets, \textit{Information and Control}, \textbf{8} (1965), 338-353.

\bibitem{Klir} G. Klir, B. Yuan, \textit{Fuzzy Sets and Fuzzy Logic}, Vol. 4, Prentice Hall, New Jersey, 1995.

\bibitem{Gomes} L. T. Gomes, L. C. de Barros, B. Bede, \textit{Fuzzy Differential Equations in Various Approaches}, Springer, Berlin, 2015.

\bibitem{Bedebook} B. Bede, \textit{Mathematics of Fuzzy Sets and Fuzzy Logic}, Springer-Verlag, Berlin, Heidelberg, 2013.

\bibitem{Leelavathi2019} R. Leelavathi, G. S. Kumar, M. S. N. Murty, Nabla Hukuhara differentiability for fuzzy functions on time scales, \textit{IAENG International Journal of Applied Mathematics}, \textbf{49} (1) (2019), 1-8.

\bibitem{LKM1} R. Leelavathi, G. Suresh Kumar, M. S. N. Murty, Second type nabla Hukuhara differentiability for fuzzy functions on time scales, \textit{Italian Journal of Pure and Applied Mathematics}, \textbf{43} (2020), 779–801.

\bibitem{Leelavathi2020} R. Leelavathi, G. Suresh Kumar, R. P. Agarwal, C. Wang, M. S. N. Murty, Generalized nabla differentiability and integrability for fuzzy functions on time scales, \textit{Axioms}, \textbf{9} (2) (2020), 65.

\bibitem{Leelavathi2018} R. Leelavathi, G. S. Kumar, M. S. N. Murty, Nabla integral for fuzzy functions on time scales, \textit{International Journal of Applied Mathematics}, \textbf{31} (5) (2018), 669.

\bibitem{Leelavathi2019-3} R. Leelavathi, G. S. Kumar, M. S. N. Murty, Characterization theorem for fuzzy functions on time scales under generalized nabla Hukuhara difference, \textit{International Journal of Innovative Technology and Exploring Engineering}, \textbf{8}(8) (2019), 1704-1706.

\bibitem{Leelavathi2019-4} R. Leelavathi, G. Suresh Kumar, M. S. N. Murty, R. V. N. Srinivasa Rao, Existence-uniqueness of solutions for fuzzy nabla initial value problems on time scales, \textit{Advances in Difference Equations}, \textbf{2019}, 1-11.

\bibitem{you} X. X. You, D. F. Zhao, B. W. Li, Nabla-Hukuhara derivative of fuzzy valued functions on time scales, \textbf{38} (4) (2018), 580--588.

\bibitem{Hari} R. Hari Kishore, R. Leelavathi, S. V. D. Mohana Rupa, A. Muneera, G. Suresh Kumar, Exploring fuzzy nabla dynamics on time scales with the characterization theorem, \textbf{28} (3) (2025), 544-553.

\bibitem{Fard} O. S. Fard, T. A. Bidgoli, Calculus of fuzzy functions on time scales (I), \textit{Soft Computing}, \textbf{19} (2015), 293-305.

\bibitem{Cano} Y. Chalco Cano, R. Rodríguez López, M. D. Jiménez Gamero, Characterizations of generalized differentiable fuzzy functions, \textit{Fuzzy Sets and Systems}, \textbf{295} (2016), 37-56.



\bibitem{qiu} D. Qiu, The generalized Hukuhara differentiability of interval valued function is not fully equivalent to the one-sided differentiability of its endpoint functions, \textit{Fuzzy Sets and Systems}, \textbf{419} (2021), 158--168.

\bibitem{longo} F. Longo, B. Laiate, M. C. Gadotti, J. F. da C. A. Meyer, Characterization results of generalized differentiabilities of fuzzy functions, \textit{Fuzzy Sets and Systems}, \textbf{490} (2024), 109038.


\bibitem{TZ} J. Tao, Z. Zhang, Properties of interval valued function space under the gH-difference and their application to semi-linear interval differential equations, \textit{Advances in Difference Equations}, (2016) 2016:45.

\bibitem{Agarwal} R. P. Agarwal, M. Bohner, Basic calculus on time scales and some of its applications, \textit{Results Math.}, \textbf{35} (1-2), 3-22, 19.


\bibitem{Negoita} C. Negoita, D. Ralescu, \textit{Application of Fuzzy Sets to System Analysis}, Wiley, New York, 1975.


\bibitem{Bede} B. Bede, L. Stefanini, Generalized differentiability of fuzzy valued functions, \textit{Fuzzy Sets and Systems}, \textbf{230} (2013), 119-141.

\bibitem{Stefanini} L. Stefanini, A generalization of Hukuhara difference and division for interval and fuzzy arithmetic, \textit{Fuzzy Sets and Systems}, \textbf{161} (11) (2010), 1564-1584.


\bibitem{Diamond2} P. Diamond, P. Kloeden, Metric spaces of fuzzy sets, \textit{Fuzzy Sets and Systems}, \textbf{35} (2) (1990), 241-249.


\bibitem{Bede1} B. Bede, J. R. Imre, A. L. Bencsik, First order linear fuzzy differential equations under generalized differentiability, \textit{Information Sciences}, \textbf{177} (2007), 1648–1662.

\end{thebibliography}
\end{document}